\documentclass[a4paper,11pt,twoside,reqno]{amsart}

\usepackage[utf8]{inputenc}
\usepackage[plainpages=false,pdfpagelabels=true]{hyperref}
\usepackage{amssymb,amsthm,wasysym}
\usepackage[margin=1in]{geometry}

\usepackage{slashed}
\usepackage[all]{xy}

\numberwithin{equation}{section}

\newtheorem{Satz}{Theorem}[section]
\newtheorem{Prop}[Satz]{Proposition}
\newtheorem{Lem}[Satz]{Lemma}

\newtheorem{Cor}[Satz]{Corollary}
\newcommand{\vol}{{\operatorname{vol}}}
\theoremstyle{definition}
\newtheorem{Dfn}[Satz]{Definition}
\newtheorem{Bem}[Satz]{Remark}
\newtheorem{Bsp}[Satz]{Example}
\newcommand{\tr}{\operatorname{tr}}
\newcommand{\hess}{\operatorname{Hess}}

\newcommand{\C}{{\mathbb{C}}}
\newcommand{\D}{\slashed{D}}

\numberwithin{equation}{section}

\title{Nonlinear Dirac equations, Monotonicity Formulas and Liouville Theorems}
\author{Volker Branding}
\date{\today}
\address{University of Vienna, Faculty of Mathematics\\
Oskar-Morgenstern-Platz 1, 1090 Vienna, Austria\\}
\email{volker.branding@univie.ac.at}
\begin{document}

\begin{abstract}
We study the qualitative behavior of nonlinear Dirac equations arising in quantum field theory
on complete Riemannian manifolds. In particular, we derive monotonicity formulas and Liouville theorems
for solutions of these equations. Finally, we extend our analysis to Dirac-harmonic maps with curvature term.
\end{abstract} 

\maketitle

\section{Introduction and results}
In quantum field theory spinors are employed to model fermions.
The equations that govern the behavior of fermions are both linear and nonlinear Dirac equations.
A Dirac equation with vanishing right hand side describes a free massless fermion
and linear Dirac equations describe free fermions having a mass.
However, to model the interaction of fermions one has to take into account nonlinearities.

In mathematical terms spinors are sections in a vector bundle, the spinor bundle, which is
defined on a Riemannian spin manifold. The spin condition is of topological nature and ensures 
the existence of the spinor bundle \(\Sigma M\).
The mathematical analysis of linear and nonlinear Dirac equations comes with two kinds of difficulties:
First of all, the Dirac operator is of first order, such that tools like the maximum principle are
not available. Secondly, in contrast to the Laplacian, the Dirac operator has its spectrum on the whole
real line.

Below we give a list of action functionals that arise in quantum field theory.
Their critical points all lead to nonlinear Dirac equations.
To this end let \(D\) be the classical Dirac operator on a Riemannian spin manifold \((M,g)\) of dimension \(n\) and \(e_i\) an orthonormal basis of \(TM\).
Furthermore, let \(\cdot\) be the Clifford multiplication of spinors with tangent vectors and \(\omega_\mathbb{C}\) the complex volume form.
Moreover, we fix a hermitian scalar product on the spinor bundle.

\begin{enumerate}
\item The \emph{Soler model} \cite{PhysRevD.1.2766} describes fermions that interact by a quartic term in the action functional.
In quantum field theory this model is usually studied on four-dimensional Minkowski space:
\begin{equation*}
E(\psi)=\int_M(\langle\psi,D\psi\rangle-\lambda|\psi|^2-\frac{\mu}{2}|\psi|^4)d\vol_g
\end{equation*}

\item The \emph{Thirring model} \cite{Thirring} describes the self-interaction of fermions in two-dimensional Minkowski space:
\begin{equation*}
E(\psi)=\int_M(\langle\psi,D\psi\rangle-\lambda|\psi|^2-\frac{\mu}{2}\sum_{i=1}^n\langle\psi,e_i\cdot\psi\rangle\langle\psi,e_i\cdot\psi\rangle)d\vol_g
\end{equation*}

\item The \emph{Nambu–Jona-Lasinio model} \cite{1961PhRv..122..345N} is a model for interacting fermions with chiral symmetry.
It also contains a quartic interaction term and is defined on an even-dimensional spacetime:
\begin{equation*}
E(\psi)=\int_M\big(\langle\psi,D\psi\rangle+\frac{\mu}{4}(|\psi|^4-\langle\psi,\omega_\C\cdot\psi\rangle\langle\psi,\omega_\C\cdot\psi\rangle)\big)d\vol_g
\end{equation*}
Note that this model does not have a term proportional to \(|\psi|^2\) in the action functional.

\item The \emph{Gross–Neveu model with N flavors} \cite{1974PhRvD..10.3235G} is a model for \(N\) interacting fermions 
in two-dimensional Minkowski space:
\begin{equation*}
E(\psi)=\int_M(\langle\psi,D\psi\rangle-\lambda|\psi|^2+\frac{\mu}{2N}|\psi|^4)d\vol_g
\end{equation*}
The spinors that we are considering here are twisted spinors, more precisely \(\psi\in\Gamma(\Sigma M\otimes\mathbb{R}^N)\).

\item The nonlinear supersymmetric sigma model in quantum field theory consists of a map \(\phi\) between two Riemannian manifolds \(M\) and \(N\)
and a spinor along that map. Moreover, \(R^N\) is the curvature tensor on \(N\) and \(\D\) 
denotes the corresponding Dirac operator. The action functional under consideration is
\begin{equation*}
E_c(\phi,\psi)=\int_M(|d\phi|^2+\langle\psi,\D\psi\rangle-\frac{1}{6}\langle R^N(\psi,\psi)\psi,\psi\rangle)d\vol_g.
\end{equation*}
The critical points of this functional became known in the mathematics literature as \emph{Dirac-harmonic maps with curvature term}.
In contrast to the physics literature this mathematical version of the nonlinear supersymmetric sigma model employs commuting
spinors while in physics anticommuting spinors are used.
\end{enumerate}
In the models (1)-(4) from above the real parameter \(\lambda\) can be interpreted as mass, whereas
the real constant \(\mu\) describes the strength of interaction.
All of the models listed above lead to nonlinear Dirac equations of the form
\begin{equation}
\label{nonlinear-dirac}
D\psi\sim\lambda\psi+\mu|\psi|^2\psi. 
\end{equation}
Note that in the physics literature Clifford multiplication is usually expressed as matrix multiplication with \(\gamma^\mu\) and the complex
volume element is referred to as \(\gamma^5\).
In contrast to the physics literature we will always assume that spinors are commuting, whereas in the physics 
literature they are mostly assumed to be Grassmann-valued.
For simplicity we will mainly focus on the Soler model.

Several existence results for equations of the form \eqref{nonlinear-dirac} are available:
In \cite{MR2813439} existence results for nonlinear Dirac equations on compact spin manifolds are obtained.
For \(n\geq 4\) existence results for nonlinear Dirac equation with critical exponent on compact spin manifolds, 
that is
\[
D\psi=\lambda\psi+|\psi|^\frac{2}{n-1}\psi
\]
with \(\lambda\in\mathbb{R}\), have been obtained in \cite{MR2733578}.
For \(\lambda=0\) this equation is known as the spinorial Yamabe equation.
In particular, this equation is interesting for \(n=2\) since it is closely related to conformally
immersed constant mean curvature surfaces in \(\mathbb{R}^3\). 
Moreover, existence results for the spinorial Yamabe equation have been obtained on \(S^3\) \cite{MR3299382}
and on \(S^n\) \cite{MR3037015} for \(n\geq 2\). For a spectral and geometric analysis of the spinorial Yamabe equation we refer to \cite{MR2550205}.
The regularity of weak solutions of equations of the form \eqref{nonlinear-dirac} can be established
with the tools from \cite{MR2661574} and \cite{MR2733578}, Appendix A.

Let us give an overview on the structure and the main results of the article:

In Section 2 we study general properties of nonlinear Dirac equations. In particular, we recall
the construction for identifying spinor bundles belonging to different metrics and use it to
derive the stress-energy tensor for the Soler model. 

In Section 3 we study nonlinear Dirac equations on closed Riemannian surfaces.
The first main result is Theorem \ref{theorem-surface-wkp} which states that for solutions of
equations of the form \eqref{nonlinear-dirac} for which the \(L^4\)-norm of \(\psi\) is 
sufficiently small on a disc \(D\) all \(W^{k,p}\) norms can be controlled on a smaller disc \(D'\subset D\).
Moreover, in Proposition \ref{proposition-nodal-set} we present an estimate on the nodal set of solutions of \eqref{nonlinear-dirac}
and Proposition \ref{proposition-surface-vanishing} shows that solutions of equations of the form \eqref{nonlinear-dirac} must be trivial 
if \(\lambda=0\) and \(|\psi|_{L^4(M)}\) is sufficiently small.

In Section 4 we investigate nonlinear Dirac equations on complete noncompact Riemannian manifolds.
First, we will prove Theorem \ref{theorem-liouville-stationary-soler} which states that stationary solutions of equations of the form \eqref{nonlinear-dirac}
with finite energy must be trivial if \(M=\mathbb{R}^n,\mathbb{H}^n,n\geq 3\).
Moreover, in Proposition \ref{monotonicity-type-rn} we show that for \(M=\mathbb{R}^n,n\geq 3\) for critical points of the Soler model
the quantity \(R^{2-n}\int_{B_R}|\psi|^4d\mu\) is almost monotone increasing in \(R\).
Moreover, we discuss the problems that arise when trying to extend the analysis to the case of
a Riemannian manifold.
Finally, in Theorem \ref{theorem-liouville-complete-soler} we show that critical points of the Soler model on a complete noncompact Riemannian
manifold with positive Ricci curvature satisfying an additional energy condition must be trivial.

In Section 5 we focus on Dirac-harmonic maps with curvature term from complete manifolds.
The latter consist of a pair of a map between two Riemannian manifolds and a vector spinor defined along that map.
First, we will show that stationary Dirac-harmonic maps with curvature term from 
\(M=\mathbb{R}^n,\mathbb{H}^n,n\geq 3\) to target spaces with positive sectional curvature must be
trivial if a certain energy is finite (Theorem \ref{theorem-dhc-stationary}).
Moreover, in the case that \(M=\mathbb{R}^n,n\geq 3\), we will establish an almost monotonicity formula 
(Proposition \ref{dhc-mono-rn}) and also discuss its extension to the case of a Riemannian manifold.
Finally, we show that Dirac-harmonic maps with curvature term from complete Riemannian manifolds
with positive Ricci curvature to target manifolds with negative sectional curvature
must be trivial if a certain energy is finite and a certain inequality relating
Ricci curvature and energy holds (Theorem \ref{theorem-liouville-dhc}).

\section{Nonlinear Dirac equations on Riemannian manifolds}
Let \((M,g)\) be a Riemannian spin manifold of dimension \(n\).
A Riemannian manifold admits a spin structure if the second Stiefel-Whitney class of its tangent bundle vanishes.

We briefly recall the basic notions from spin geometry,
for a detailed introduction to spin geometry we refer to the book \cite{MR1031992}.

We fix a spin structure on the manifold \(M\) and consider the spinor bundle \(\Sigma M\).
On the spinor bundle \(\Sigma M\) we have the Clifford multiplication of spinors with tangent vectors denoted by \(\cdot\).
Moreover, we fix a hermitian scalar product on the spinor bundle and denote its real part by \(\langle\cdot,\cdot\rangle\).
Clifford multiplication is skew-symmetric 
\[
\langle\psi,X\cdot\xi\rangle=-\langle X\cdot\psi,\xi\rangle
\]
for all \(\psi,\xi\in\Gamma(\Sigma M)\) and \(X\in TM\).
Moreover, the Clifford relations
\begin{equation}
\label{clifford}
X\cdot Y+Y\cdot X=-2g(X,Y)
\end{equation}
hold for all \(X,Y\in TM\).
The Dirac operator \(D\colon\Gamma(\Sigma M)\to\Gamma(\Sigma M)\)
is defined as the composition of first applying the covariant derivative on the spinor bundle
followed by Clifford multiplication. More precisely, it is given by
\[
D:=\sum_{i=1}^ne_i\cdot\nabla^{\Sigma M}_{e_i},
\]
where \(e_i,i=1\ldots n\) is an orthonormal basis of \(TM\).
Sometimes we will make use of the Einstein summation convention and just sum over repeated indices.
The Dirac operator is of first order, elliptic and self-adjoint with respect to the \(L^2\)-norm. 
Hence, if \(M\) is compact the Dirac operator has a real and discrete spectrum.

The square of the Dirac operator satisfies the \emph{Schroedinger-Lichnerowicz} formula
\begin{equation}
\label{schroedinger}
D^2=\nabla^\ast\nabla+\frac{R}{4},
\end{equation}
where \(R\) denotes the scalar curvature of the manifold \(M\).

After having recalled the basic definitions from spin geometry we will focus on the analysis of the following
action functional (which is the first one from the introduction)
\begin{equation}
\label{soler-energy}
E(\psi)=\int_M(\langle\psi,D\psi\rangle-\lambda|\psi|^2-\frac{\mu}{2}|\psi|^4)d\vol_g.
\end{equation}
Its critical points are given by
\begin{equation}
\label{soler}
D\psi=\lambda\psi+\mu|\psi|^2\psi.
\end{equation}
It turns out that \(L^4(\Sigma M)\times W^{1,\frac{4}{3}}(\Sigma M)\) is the right function space
for weak solutions of \eqref{soler}.

\begin{Dfn}
We call \(\psi\in L^4(\Sigma M)\times W^{1,\frac{4}{3}}(\Sigma M)\) a \emph{weak solution}
if it solves \eqref{soler} in a distributional sense.
\end{Dfn}

The analytic structure of the other action functionals listed in the introduction is the same as 
the one of \eqref{soler-energy}. Due to this reason many of the results that will be obtained
for solutions of \eqref{soler} can easily be generalized to critical points of the other models.

The equation \eqref{soler} is also interesting from a geometric point of view since it interpolates
between eigenspinors (\(\mu=0\)) and a non-linear Dirac equation (\(\lambda=0\)) that arises
in the study of CMC immersions from surfaces into \(\mathbb{R}^3\).

In the following we want to vary the action functional \eqref{soler-energy} 
(and later on also other similar functionals) with respect to the metric \(g\).
There had been many isolated mathematical results in the literature how to carry out this calculation before a
first complete framework for the Riemannian case was established in \cite{MR1158762}.
Later, this was generalized to the pseudo-Riemannian case in \cite{MR2121740}.

We will now give a brief survey on how to identify spinor bundles belonging to different metrics
recalling the methods that were established in \cite{MR1158762}.
However, our presentation of these methods is motivated from the one of \cite{MR1738150}, Chapter 2.

Suppose we have two spinor bundles \(\Sigma_g M\) and \(\Sigma_hM\) corresponding
to different metrics \(g\) and \(h\). There exists a unique positive definite tensor field
\(h_{g}\) uniquely determined by the requirement \(h(X,Y)=g(HX,HY)=g(X,h_{g}Y)\),
where \(H:=\sqrt{h_g}\).
Let \(P_g\) and \(P_h\) be the oriented orthonormal frame bundles of \((M,g)\) and \((M,h)\).
Then \(H^{-1}\) induces an equivariant isomorphism \(b_{g,h}\colon P_g\to P_{h}\) via the assignment
\(E_i\mapsto H^{-1}E_i,i=1\ldots n\).
We fix a spin structure \(\Lambda_g\colon Q_g\to P_g\) of \((M,g)\) and
think of it as a \(\mathbb{Z}_2\)-bundle. The pull-back of \(\Lambda_g\)
via the isomorphism \(b_{h,g}\colon P_h\to P_{g}\) induces a \(\mathbb{Z}_2\)-bundle
\(\Lambda_h\colon Q_h\to P_h\). Moreover, we get a Spin\((n)\)-equivariant isomorphism 
\(\tilde{b}_{h,g}\colon Q_h\to Q_g\) such that the following diagram commutes:

$$
\xymatrix{ 
Q_h \ar[d]_{\Lambda_h} \ar[r]^{\tilde{b}_{h,g}} & Q_{g} \ar[d]^{\Lambda_{g}} &\\
P_h\ar[r]^{b_{h,g}} & P_{g}&}
$$
Making use of this construction we obtain the following
\begin{Lem}
There exist natural isomorphisms 
\[
b_{g,h}\colon TM\to TM,\qquad \beta_{g,h}\colon \Sigma_gM\to\Sigma_{h}M
\]
that satisfy
\begin{align*}
h(b_{h,g}X,b_{h,g}Y)=&g(X,Y), \qquad  
\langle\beta_{h,g}\chi,\beta_{h,g}\psi\rangle_{\Sigma_hM}=\langle\psi,\chi\rangle_{\Sigma_{g}M}, \\
(b_{g,h}X)\cdot(\beta_{g,h}\psi)=&\beta_{g,h}(X\cdot\psi)
\end{align*}
for all \(X,Y\in\Gamma(TM)\) and \(\psi,\chi\in\Gamma(\Sigma_{g}M)\).
\end{Lem}

In order to calculate the variation of the Dirac operator with respect to the metric we need the following objects:
Let \(\operatorname{Sym}(0,2)\) be the space of all symmetric \((0,2)\)-tensor fields on \((M,g)\).
Any element \(k\) of \(\operatorname{Sym}(0,2)\) induces a \((1,1)\)-tensor field \(k_g\) via
\(k(X,Y)=g(X,k_gY)\). We denote the Dirac operator on \((M,g+tk)\) by \(D_{g+tk}\) for a small parameter \(t\).
Moreover, we will use the notation \(\psi_{g+tk}:=\beta_{g,g+tk}\psi\in\Gamma(\Sigma M_{g+tk})\),
which can be thought of as push-forward of \(\psi\in\Sigma_gM\) to \(\psi\in\Sigma_{g+tk}M\).
Applying the technical construction outlined above let us now recall the following classic result
from \cite{MR1158762}:

\begin{Lem}
The variation of the Dirac-energy with respect to the metric is given by
\begin{equation}
\label{variation-dirac}
\frac{d}{dt}\big|_{t=0}\langle\psi_{g+tk},D_{g+tk}\psi_{g+tk}\rangle_{\Sigma_{g+tk}M}
=-\frac{1}{4}\langle e_i\cdot\nabla^{\Sigma M}_{e_j}\psi+e_j\cdot\nabla^{\Sigma M}_{e_i}\psi,\psi\rangle_{\Sigma_gM}k^{ij},
\end{equation}
where the tensor on the right hand side is the stress-energy tensor associated to the Dirac energy.
\end{Lem}
\begin{proof}
A proof can be found in \cite[Section III]{MR1158762}.
\end{proof}

\begin{Dfn}
A weak solution \(\psi\in L^4(\Sigma M)\times W^{1,\frac{4}{3}}(\Sigma M)\) of \eqref{soler} is called \emph{stationary}
if it is also a critical point of \(E(\psi)\) with respect
to domain variations.
\end{Dfn}

\begin{Prop}
A stationary solution \(\psi\in L^4(\Sigma M)\times W^{1,\frac{4}{3}}(\Sigma M)\) of \eqref{soler-energy} satisfies
\begin{equation}
\label{stationary-identity}
0=\int_M(\langle e_i\cdot\nabla^{\Sigma M}_{e_j}\psi+e_j\cdot\nabla^{\Sigma M}_{e_i}\psi,\psi\rangle-g_{ij}\mu|\psi|^4)k^{ij}d\vol_g,
\end{equation}
where \(k^{ij}\) is a smooth element of \(\operatorname{Sym}(0,2)\).
\end{Prop}
\begin{proof}
Let \(k\) be a symmetric \((0,2)\)-tensor and \(t\) some small number.
Recall that the variation of the volume-element is given by
\begin{equation}
\label{variation-volume}
\frac{d}{dt}\big|_{t=0}d\vol_{g+tk}=\frac{1}{2}\langle g,k\rangle_gd\vol_{g}. 
\end{equation}
Moreover, as \(\beta_{g,g+tk}\) acts as an isometry on the spinor bundle we obtain
\[
|\psi_{g+tk}|_{\Sigma_{g+tk}M}^2=|\beta_{g,g+tk}\psi|_{\Sigma_{g+tk}M}^2=|\psi|_{\Sigma_gM}^2.
\]
Now, we calculate
\begin{align*}
\frac{d}{dt}\big|_{t=0}&\int_M\big(\langle\psi_{g+tk},D_{g+tk}\psi_{g+tk}\rangle_{\Sigma_{g+tk}M}
-\lambda|\psi_{g+tk}|_{\Sigma_{g+tk}M}^2-\frac{\mu}{2}|\psi_{g+tk}|_{\Sigma_{g+tk}M}^4\big)d\vol_{g+tk}\\
=&-\frac{1}{4}\int_M\langle e_i\cdot\nabla^{\Sigma M}_{e_j}\psi+e_j\cdot\nabla^{\Sigma M}_{e_i}\psi,\psi\rangle k^{ij}d\vol_g\\
&+\frac{1}{2}\int_M(\langle\psi,D\psi\rangle-\lambda|\psi|^2-\frac{\mu}{2}|\psi|^4)\langle g,k\rangle_g d\vol_g \\
=&-\frac{1}{4}\int_M\langle e_i\cdot\nabla^{\Sigma M}_{e_j}\psi+e_j\cdot\nabla^{\Sigma M}_{e_i}\psi,\psi\rangle k^{ij}d\vol_g
+\frac{1}{4}\int_M\mu|\psi|^4\langle g,k\rangle_g d\vol_g,
\end{align*}
where we used \eqref{variation-dirac} in the first step and the equation for the spinor \(\psi\), that is \eqref{soler},
in the second step completing the proof.
\end{proof}

For a smooth solution \(\psi\) of \eqref{soler} we thus obtain the stress-energy tensor
\begin{equation}
\label{stress-energy-soler}
S_{ij}=\langle e_i\cdot\nabla^{\Sigma M}_{e_j}\psi+e_j\cdot\nabla^{\Sigma M}_{e_i}\psi,\psi\rangle-g_{ij}\mu|\psi|^4.
\end{equation}

Its trace can easily be computed to be
\[
\tr S=g^{ij}S_{ij}=2\lambda|\psi|^2+(2-n)\mu|\psi|^4.
\]
Note that the stress-energy tensor is traceless for \(\lambda=0\) and \(n=2\) since
it arises from a conformally invariant action functional in that case.

\begin{Lem}
\label{lemma-stress-energy-soler}
Suppose that \(\psi\) is a smooth solution of \eqref{soler}.
Then the stress-energy tensor \eqref{stress-energy-soler} is symmetric and divergence-free.
\end{Lem}
\begin{proof}
We choose a local orthonormal basis of \(TM\) such that \(\nabla_{e_i}e_j=0, i,j=1,\ldots,n\) at the considered point.
To show that the stress-energy tensor is divergence-free we calculate
\begin{align*}
\nabla^jS_{ij}=&\nabla^j(\langle e_i\cdot\nabla_{e_j}\psi+e_j\cdot\nabla_{e_i}\psi,\psi\rangle-g_{ij}\mu|\psi|^4) \\
=&\langle e_i\cdot\Delta\psi,\psi\rangle+
\underbrace{\langle e_i\cdot\nabla_{e_j}\psi,\nabla_{e_j}\psi	\rangle}_{=0}+
\langle D\nabla_{e_i}\psi,\psi\rangle-\langle\nabla_{e_i}\psi,D\psi\rangle
-4\mu|\psi|^2\langle\nabla_{e_i}\psi,\psi\rangle.
\end{align*}
By a direct computation we find
\begin{align*}
\langle D\nabla_{e_i}\psi,\psi\rangle=&\underbrace{\langle\psi,e_j\cdot R^{\Sigma M}(e_j,e_i)\psi\rangle}_{=\frac{1}{2}\langle\psi,\operatorname{Ric}(e_i)\cdot\psi\rangle=0} 
+\langle\nabla_{e_i}D\psi,\psi\rangle 
=(\lambda+3\mu|\psi|^2)\langle\nabla_{e_i}\psi,\psi\rangle, \\
\langle\nabla_{e_i}\psi,D\psi\rangle=&(\lambda+\mu|\psi|^2)\langle\nabla_{e_i}\psi,\psi\rangle,
\end{align*}
where we used that \(\psi\) is a solution of \eqref{soler}. Thus, we obtain
\[
\nabla^jS_{ij}=\langle e_i\cdot\Delta\psi,\psi\rangle-2\mu|\psi|^2\langle\nabla_{e_i}\psi,\psi\rangle.
\]
Using \eqref{schroedinger} and \eqref{soler} we find that
\[
\langle e_i\cdot\Delta\psi,\psi\rangle=-\mu\langle e_i\cdot(\nabla|\psi|^2)\cdot\psi,\psi\rangle=\mu g(e_i,\nabla|\psi|^2)|\psi|^2=2\mu|\psi|^2\langle\nabla_{e_i}\psi,\psi\rangle,
\]
which completes the proof.
\end{proof}

\begin{Bem}
Every smooth solution of \eqref{soler} is also stationary.
We will give a short proof of this statement where we reverse the calculation performed
in the proof of Lemma \ref{lemma-stress-energy-soler}.

Hence, suppose we have a smooth solution of \eqref{soler}.
Differentiating \eqref{soler} with respect to \(e_i\) and taking the scalar product with \(\psi\) we find
\begin{align*}
0=&\langle\psi,\nabla_{e_i}D\psi\rangle-\lambda\langle \psi,\nabla_{e_i}\psi\rangle-\mu(\nabla_{e_i}|\psi|^2)|\psi|^2-\mu|\psi|^2\langle\nabla_{e_i}\psi,\psi\rangle\\
=&\langle\psi,\nabla_{e_i}D\psi\rangle-\langle D\psi,\nabla_{e_i}\psi\rangle-\mu(\nabla_{e_i}|\psi|^2)|\psi|^2.
\end{align*}
Recall that for a solution of \eqref{soler} we have \(\langle e_i\cdot\Delta\psi,\psi\rangle=\frac{\mu}{2}\nabla_{e_i}|\psi|^4\) 
and together with the arguments used in the proof of Lemma \ref{lemma-stress-energy-soler} this leads to
\begin{align*}
0=\langle\psi,D\nabla_{e_i}\psi\rangle+\langle\nabla_{e_j}\psi,e_j\cdot\nabla_{e_i}\psi\rangle+\langle e_i\cdot\Delta\psi,\psi\rangle-\mu\nabla_{e_i}|\psi|^4=\nabla^jS_{ij}.
\end{align*}
Testing this equation with a smooth function \(Y\) and using integration by parts we obtain
\begin{align*}
0=\int_M\nabla^jY^iS_{ij}d\vol_{g},
\end{align*}
which is exactly the condition of being stationary \eqref{stationary-identity}.
\end{Bem}

We will often make use of the following Bochner-type equation
\begin{Lem}
Let \(\psi\) be a smooth solution of \eqref{soler}. Then the following formula holds
\begin{equation}
\label{bochner}
\Delta\frac{1}{2}|\psi|^4=\big|d|\psi|^2\big|^2+|\psi|^2|\nabla\psi|^2
+|\psi|^4\big(\frac{R}{4}-(\lambda+\mu|\psi|^2)^2\big).
\end{equation}
\end{Lem}
\begin{proof}
By a direct calculation we find
\begin{align*}
\Delta\frac{1}{2}|\psi|^4=&\big|d|\psi|^2\big|^2+|\psi|^2|\nabla\psi|^2
+\frac{R}{4}|\psi|^4-|\psi|^2\langle\psi,D^2\psi\rangle,
\end{align*}
where we used \eqref{schroedinger}. Moreover, we obtain
\begin{align}
\label{bochner-d-quadrat}
\nonumber\langle\psi,D^2\psi\rangle=&\langle\psi,D(\lambda\psi)\rangle+\langle\psi,D(\mu|\psi|^2\psi)\rangle \\
\nonumber=&\lambda^2|\psi|^2+\lambda\mu|\psi|^4
+\mu\underbrace{\langle\psi,(\nabla|\psi|^2)\cdot\psi\rangle}_{=0}+\mu\lambda|\psi|^4+\mu^2|\psi|^6\\
=&|\psi|^2(\lambda+\mu|\psi|^2)^2,
\end{align}
where we used that \(\psi\) is a solution of \eqref{soler}.
\end{proof}
Let us recall the following definitions:
\begin{Dfn}
A spinor \(\psi\in\Gamma(\Sigma M)\) is called \emph{twistor spinor} if it satisfies
\begin{equation}
\label{twistor}
\nabla^{\Sigma M}_X\psi+\frac{1}{n}X\cdot D\psi=0
\end{equation}
for all vector fields \(X\). The spinor \(\psi\) is called \emph{Killing spinor} if
it is both a twistor spinor and an eigenspinor of the Dirac operator, that is
\begin{equation}
\label{killing}
\nabla^{\Sigma M}_X\psi+\alpha X\cdot\psi=0
\end{equation}
with \(\alpha\in\mathbb{R}\).
\end{Dfn}

It is well known that Killing spinors have constant norm, that is \(|\psi|^2=const\).
However, here we have the following 
\begin{Lem}
Suppose that \(\psi\) is a solution of \eqref{soler} and
a twistor spinor. Then \(\psi\) has constant norm.
\end{Lem}
\begin{proof}
We calculate for an arbitrary \(X\in TM\)
\begin{align*}
\partial_X\frac{1}{2}|\psi|^2=\langle\nabla^{\Sigma M}_X\psi,\psi\rangle=-\frac{1}{n}\langle X\cdot D\psi,\psi\rangle
=-\frac{1}{n}(\lambda+\mu|\psi|^2)\langle X\cdot\psi,\psi\rangle,
\end{align*}
where we first used that \(\psi\) is a twistor spinor and then used that \(\psi\) is a solution of \eqref{soler}. 
The statement then follows from the skew-symmetry of the Clifford multiplication.
\end{proof}

\begin{Bsp}
Suppose that \(\psi\) is a Killing spinor with constant \(\alpha=\frac{\lambda+\mu|\psi|^2}{n}\).
Then it is a solution of \eqref{soler}.
However, this above approach is rather restrictive since only few Riemannian manifolds
admit Killing spinors \cite{MR1224089}.
\end{Bsp}

\begin{Prop}
Suppose that \(\psi\) is a smooth solution of \eqref{soler} and also a twistor spinor. 
Then the stress-energy tensor \eqref{stress-energy-soler} acquires the form
\begin{equation}
S_{ij}=\frac{1}{n}g_{ij}\big(\mu(2-n)|\psi|^4+2\lambda|\psi|^2\big).
\end{equation}
\end{Prop}
In particular, the stress-energy tensor is just a multiple of the metric.
\begin{proof}
We consider the stress-energy tensor \eqref{stress-energy-soler}
and use the fact that \(\psi\) is a twistor spinor, that is
\begin{align*}
S_{ij}=&\langle e_i\cdot\nabla_{e_j}\psi+e_j\cdot\nabla_{e_i}\psi,\psi\rangle-g_{ij}\mu|\psi|^4 \\
=&-\frac{1}{n}\langle (\underbrace{e_i\cdot e_j+e_j\cdot e_i}_{=-2g_{ij}})D\psi,\psi\rangle-g_{ij}\mu|\psi|^4
=\mu(\frac{2}{n}-1)|\psi|^4g_{ij}+\frac{2}{n}\lambda|\psi|^2 g_{ij},
\end{align*}
where we used the Clifford relations \eqref{clifford} and \eqref{soler}.
\end{proof}

\section{Nonlinear Dirac equations on closed surfaces}
In this section we will derive several properties of solutions of \eqref{soler} on closed Riemannian surfaces.
First, we derive a local energy estimate for smooth solutions of \eqref{soler}.
Our result is similar to the energy estimate that was obtained in \cite{MR2390834}, Theorem 2.1,
which corresponds to \eqref{soler} with \(\lambda=0\).
We obtain the following
\begin{Satz}
\label{theorem-surface-wkp}
Let \(\psi\) be a smooth solution of \eqref{soler}. If \(|\psi|_{L^4(D)}<\epsilon\) then
\begin{equation}
|\psi|_{W^{k,p}(D')}\leq C|\psi|_{L^4(D)}
\end{equation}
for all \(D'\subset D\) and \(p>1\). The constant \(C\) depends on \(D',\mu,\lambda,k,p\).
\end{Satz}

The statement of the above Theorem would also hold true if \(\psi\) was only a weak solution of \eqref{soler},
that is \(\psi\in L^4(\Sigma M)\times W^{1,\frac{4}{3}}(\Sigma M)\).
By the regularity theory presented in \cite{MR2661574} a distributional solution of \eqref{soler}
with \(\psi\in L^4(\Sigma M)\) is actually smooth if \(\dim M=2\).

We will divide the proof into two Lemmas, the result then follows by iterating the procedure
outlined below.

\begin{Lem}
Let \(\psi\) be a smooth solution of \eqref{soler}. If \(|\psi|_{L^{4}}(D)<\epsilon\) then
for all \(p>1\) and all \(D'\subset D\) we have
\begin{align}
\label{psi-local-l4}
|\psi|_{L^{p}(D')}\leq C|\psi|_{L^{4}(D)},
\end{align}
where the constant \(C\) depends on \(D',\mu,\lambda,p\).
\end{Lem}
\begin{proof}
Choose a cut-off function \(\eta\) with \(0\leq\eta\leq 1\), \(\eta|_{D'}=1\) and \(\operatorname{supp}\eta\subset D\).
Then we have
\begin{align*}
D(\eta\psi)=&\eta D\psi+\nabla\eta\cdot\psi 
=\eta\lambda\psi+\eta\mu|\psi|^2\psi+\nabla\eta\cdot\psi.
\end{align*}
We set \(\xi=\eta\psi\) and by making use of elliptic estimates for first order equations we obtain
\begin{align*}
|\xi|_{W^{1,q}(D)}\leq &C(|\eta\psi|_{L^q(D)}+\mu\big||\psi^3\eta|\big|_{L^q(D)}+|\nabla\eta||\psi|_{L^q(D)}) \\
\leq& C(|\psi|_{L^q(D)}+\big||\psi^3\eta|\big|_{L^q(D)}).
\end{align*}
We set \(q^\ast:=\frac{2q}{2-q}\) for \(q<2\). By the Hölder inequality we get
\[
\big||\psi^3\eta|\big|_{L^q(D)}\leq|\psi|^2_{L^{4}(D)}|\xi|_{L^{q^\ast}(D)}.
\]
Applying the Sobolev embedding theorem in two dimensions we find
\[
|\xi|_{L^{q^\ast}(D)}\leq C|\xi|_{W^{1,q}(D)}\leq C(|\psi|_{L^q(D)}+|\psi|^2_{L^{4}(D)}|\xi|_{L^{q^\ast}(D)}).
\]
Using the small energy assumption we get
\[
|\xi|_{L^{q^\ast}(D)}\leq C|\psi|_{L^{4}(D)}.
\]
For any \(p>1\) we can find some \(q<2\) such that \(p=q^\ast\).
\end{proof}

\begin{Lem}
Let \(\psi\) be a smooth solution of \eqref{soler}. If \(|\psi|_{L^{4}}(D)<\epsilon\) then
for all \(p>1\) and all \(D'\subset D\) we have
\begin{align}
|\psi|_{W^{1,p}(D')}\leq C|\psi|_{L^4(D)},
\end{align}
where the constant \(C\) depends on \(D',\mu,\lambda,p\).
\end{Lem}
\begin{proof}
Again, choose a cut-off function \(\eta\) with \(0\leq\eta\leq 1\), \(\eta|_{D'}=1\) and \(\operatorname{supp}\eta\subset D\).
Setting \(\xi=\eta\psi\) we locally have
\[
\int_D|\nabla\xi|^2dx=\int_D|D\xi|^2dx=\int_D|\eta\lambda\psi+\eta\mu|\psi|^2\psi+\nabla\eta\cdot\psi|^2dx\leq C\int_D(|\psi|^2+|\psi|^6)dx.
\]
We obtain the following inequality
\[
|\nabla\xi|_{L^2(D)}\leq C(|\psi|^3_{L^6(D)}+|\psi|_{L^2(D)})\leq C|\psi|_{L^4(D)},
\]
which yields
\begin{equation}
\label{psi-local-nabla-psi-l2}
|\nabla\psi|_{L^2(D')}\leq C|\psi|_{L^4(D)}.
\end{equation}
By a direct computation we find
\[
D^2\psi=\lambda^2\psi+2\mu\lambda|\psi|^2\psi+\mu^2|\psi|^4\psi+\mu(\nabla|\psi|^2)\cdot\psi
\]
and also
\[
\Delta\xi=(\Delta\eta)\psi+2\nabla\eta\nabla\psi+\eta\Delta\psi.
\]
This yields
\[
|\Delta\xi|\leq C(|\psi|+|\nabla\psi|+|\psi|^2|\nabla\psi|+|\psi|^3+|\psi|^5).
\]
On the disc \(D\) we have \(\Delta=-D^2\), hence we find
\begin{equation}
\label{psi-local-w2p}
|\eta\psi|_{W^{2,p}(D)}\leq C(|\psi|_{L^p(D)}+|\nabla\psi|_{L^p(D)}+\big||\psi|^2|\nabla\psi|\big|_{L^p(D)}
+\big||\psi|^3\big|_{L^p(D)}+\big||\psi|^5\big|_{L^p(D)}).
\end{equation}
Using \eqref{psi-local-l4} and \eqref{psi-local-nabla-psi-l2} we obtain
\[
\big||\psi|^2|\nabla\psi|\big|_{L^p(D)}\leq C|\nabla\psi|_{L^2(D')}|\psi|^2_{L^8(D')}\leq C|\psi|_{L^4(D)}
\]
and the same bound applies to the first and the last two terms of \eqref{psi-local-w2p}.
Thus, we obtain by setting \(p=\frac{4}{3}\) in \eqref{psi-local-l4} and applying the Sobolev embedding theorem
\[
|\psi|_{W^{1,4}(D')}\leq C|\eta\psi|_{W^{2,\frac{4}{3}}(D')}\leq C|\psi|_{L^4(D)}
\]
for all \(D'\subset D\).
In particular, this implies
\[
|\psi|_{L^\infty(D')}\leq C|\psi|_{L^4(D)}.
\]
At this point we may set \(p=2\) in \eqref{psi-local-w2p} and find
\[
|\psi|_{W^{1,p}(D')}\leq C|\psi|_{W^{2,2}(D')}\leq C|\psi|_{L^4(D)},
\]
which proves the result.
\end{proof}

\begin{Bem}
In the case that \(\lambda=0\) the equation \eqref{soler} arises from a conformally invariant action functional and
is scale invariant. This scale invariance can be exploited to show that solutions of \eqref{soler}
cannot have isolated singularities, see \cite{MR2390834}, Theorem 3.1.
\end{Bem}

By the main result of \cite{MR1714341} we know that the nodal set of solutions to \eqref{soler} on closed surfaces is discrete.
The next Proposition gives an upper bound on their nodal set.
\begin{Prop}
\label{proposition-nodal-set}
Suppose that \(\psi\) is a smooth solution of \eqref{soler} that is not identically zero. 
Then the following inequality holds
\begin{equation}
\label{nodal-set}
\int_M(\lambda+\mu|\psi|^2)^2d\vol_g\geq 2\pi\chi(M)+4\pi N,
\end{equation}
where \(\chi(M)\) is the Euler characteristic of the surface. Moreover,
\(N\) denotes an estimate on the nodal set
\[
N=\sum_{p\in M,|\psi|(p)=0}n_p, 
\]
where \(n_p\) is the order of vanishing of \(|\psi|\) at the point \(p\).
\end{Prop}
\begin{proof}
Throughout the proof we assume that \(\psi\neq 0\).
Now, we recall the following inequality (see \cite{MR3774399}, Lemma 2.1 and references therein for a detailed derivation)
\[
\frac{\langle\psi,D^2\psi\rangle}{|\psi|^2}\geq\frac{R}{4}+\frac{|T|^2}{4|\psi|^4}-\Delta\log|\psi|
\]
with the stress-energy tensor for the Dirac action \(T(X,Y):=\langle X\cdot\nabla_Y\psi+Y\cdot\nabla_X\psi,\psi\rangle\).
Using \eqref{bochner-d-quadrat} we find
\[
\frac{\langle\psi,D^2\psi\rangle}{|\psi|^2}=(\lambda+\mu|\psi|^2)^2.
\]
We can estimate the stress-energy tensor as
\[
|T|^2\geq 2(\lambda+\mu|\psi|^2)^2,
\]
which yields
\[
(\lambda+\mu|\psi|^2)^2\geq K-2\Delta\log|\psi|,
\]
where \(K=2R\) denotes the Gaussian curvature of \(M\).
By integrating over \(M\) and using that for a function with discrete zero set 
\[
\int_M\Delta\log|\psi|d\vol_g=-2\pi\sum_{p\in M,|\psi|(p)=0}n_p
\]
we obtain the result.
\end{proof}

\begin{Bem}
The estimate on the nodal set \eqref{nodal-set} generalizes the estimates on the nodal set for eigenspinors \cite{MR3774399} 
and on solutions to non-linear Dirac equations \cite{MR2550205}, Proposition 8.4.
\end{Bem}

\begin{Cor}
\begin{enumerate}
 \item Due to the last Proposition we obtain the following upper bound on the nodal set of solutions to \eqref{soler}
 \[
  N\leq-\frac{\chi(M)}{2}+\frac{1}{4\pi}\int_M(\lambda+\mu|\psi|^2)^2d\vol_g.
 \]
 \item We also obtain a vanishing result for surfaces of positive Euler characteristic:
More precisely, if
\begin{equation*}
\int_M(\lambda+\mu|\psi|^2)^2d\vol_g<4\pi
\end{equation*}
then we get a contradiction from \eqref{nodal-set} forcing \(\psi\) to be trivial.
\end{enumerate}
\end{Cor}

Using the Sobolev embedding theorem we can obtain another variant of the last statement from the previous Corollary.
\begin{Prop}
\label{proposition-surface-vanishing}
Let \(\psi\) be a smooth solution of \eqref{soler} with \(\lambda=0\).
Suppose that there do not exist harmonic spinors on \(M\). 
There exists some \(\epsilon_0>0\) depending on \(M,\mu\) such that whenever \(\epsilon<\epsilon_0\) and
\begin{equation}
|\psi|^2_{L^4}<\epsilon
\end{equation}
we have \(\psi=0\).
\end{Prop}
\begin{proof}
By assumption \(0\) is not in the spectrum of \(D\) and we can estimate
\[
|\psi|\leq\frac{1}{|\lambda_1|}|D\psi|,
\]
where \(\lambda_1\) denotes the smallest eigenvalue of the Dirac operator.
Making use of elliptic estimates for first order equations we find
\begin{align*}
|\psi|_{L^4}\leq C|\psi|_{W^{1,\frac{4}{3}}} 
&\leq C(|D\psi|_{L^\frac{4}{3}}+|\psi|_{L^\frac{4}{3}})\\
&\leq C|\mu|\big||\psi|^3\big|_{L^\frac{4}{3}}\\
&\leq C|\mu||\psi|^3_{L^4}\\
&\leq \epsilon C|\mu||\psi|_{L^4},
\end{align*}
where we made use of the assumptions. Thus, for \(\epsilon\) small enough \(\psi\) has to vanish.
\end{proof}

\begin{Bem}
The regularity theory for Dirac-type equations on Riemannian manifolds is well-established,
see for example the \(L^2\)-theory developed in \cite{MR2124018}. 
Recently, it could be substantially extended in \cite{MR3908762} to also include higher \(L^p\)-norms.
Using this recent regularity theory for Dirac equations \cite[Theorem 1.1]{MR3908762}
it should be possible to get rid of the requirement that \(M\) is not supposed to admit harmonic spinors
in Proposition \ref{proposition-surface-vanishing}. However, Theorem 1.1 in \cite{MR3908762}
is formulated for boundary value problems of Dirac-type operators and it would be necessary
to obtain a variant of this result for closed manifolds. 
Having such a result at hand the proof of Proposition \ref{proposition-surface-vanishing}
could be simplified in such a way that one does not need the condition of \(M\) having
no harmonic spinors.

However, it seems that a variant of \cite[Theorem 1.1]{MR3908762} on closed manifolds, which would be a global statement,
could not help to improve Theorem \ref{theorem-surface-wkp} as this theorem is of a local nature and Proposition \ref{proposition-surface-vanishing}
shows that demanding \(|\psi|_{L^4}<\epsilon\) globally forces \(\psi\) to be trivial.
\end{Bem}

\subsection{The higher-dimensional case}
\begin{Prop}
Suppose that \(M\) is a closed Riemannian spin manifold with positive scalar curvature.
Suppose that \(\psi\) is a smooth solution of \eqref{soler} with small energy,
that is
\begin{equation}
(\lambda+\mu|\psi|^2)^2<\frac{R}{4}.
\end{equation}
Then \(\psi\) vanishes identically.
\end{Prop}
\begin{proof}
We use the Bochner formula \eqref{bochner} and calculate
\[
\Delta\frac{1}{2}|\psi|^4=\big|d|\psi|^2\big|^2+|\psi|^2|\nabla\psi|^2
+|\psi|^4\big(\frac{R}{4}-(\lambda+\mu|\psi|^2)^2\big)>0
\]
using the assumption. Hence \(|\psi|^4\) is a subharmonic function and due to the maximum principle 
it has to be constant.
Thus, we obtain
\[
0=|\psi|^2|\nabla\psi|^2+|\psi|^4\big(\frac{R}{4}-(\lambda+\mu|\psi|^2)^2\big)
\]
and the result follows by making use of the assumption.
\end{proof}

\section{Nonlinear Dirac equations on complete manifolds}
In this section we study the behavior of solutions of \eqref{soler} on complete manifolds.
We will derive several monotonicity formulas and, as an application, we obtain 
Liouville theorems.
\subsection{A Liouville Theorem for stationary solutions}
In this section we will derive a vanishing theorem for stationary solutions of \eqref{soler}.
\begin{Satz}
\label{theorem-liouville-stationary-soler}
Suppose that \(M=\mathbb{R}^n, \mathbb{H}^n\) with \(n\geq 3\).
Let \(\psi\in L_{loc}^4(\Sigma M)\times W_{loc}^{1,\frac{4}{3}}(\Sigma M)\) be a stationary solution of \eqref{soler}.
If \(\lambda\mu\leq 0\) and 
\begin{equation}
\int_M(|\psi|^4+|\nabla\psi|^\frac{4}{3})d\vol_g<\infty
\end{equation}
then \(\psi\) vanishes identically.
\end{Satz}
\begin{proof}
We will first show the result for \(M=\mathbb{R}^n\).
Choose \(\eta\in C^\infty_0(\mathbb{R})\) such that \(\eta=1\) for \(r\leq R\), 
\(\eta=0\) for \(r\geq 2R\) and \(|\eta'(r)|\leq\frac{4}{R}\). In addition, we choose \(Y(x)=x\eta(r)\in C^\infty(M,\mathbb{R}^n)\),
where \(r=|x|\).
Then, we set
\[
k_{ij}:=\frac{\partial Y_i}{\partial x^j}=\eta(r)\delta_{ij}+\frac{x_ix_j}{r}\eta'(r)
\]
and inserting this into \eqref{stationary-identity} we obtain
\begin{align*}
\int_{\mathbb{R}^n}(2\langle\psi,D\psi\rangle-n\mu|\psi|^4)\eta(r)d\vol_g
=-\int_{\mathbb{R}^n}(2\langle\psi,\partial_r\cdot\nabla_{\partial_r}\psi\rangle-\mu|\psi|^4)r\eta'(r)d\vol_g.
\end{align*}
Using the equation for \(\psi\) we get
\begin{align*}
\int_{\mathbb{R}^n}(2\lambda|\psi|^2+\big(2-n)\mu|\psi|^4)\eta(r)d\vol_g
=-\int_{\mathbb{R}^n}(2\langle\psi,\partial_r\cdot\nabla_{\partial_r}\psi\rangle-\mu|\psi|^4)r\eta'(r)d\vol_g.
\end{align*}
The right hand side can be controlled as follows
\[
\big|\int_{\mathbb{R}^n}(2\langle\psi,\partial_r\cdot\nabla_{\partial_r}\psi\rangle-\mu|\psi|^4)r\eta'(r)d\vol_g\big|\leq
C\int_{B_{2R}\setminus B_R}(|\psi||\nabla\psi|+|\psi|^4)dx.
\]
First, we consider the case that \(\lambda\geq 0\) and \(\mu\leq 0\).
Making use of the assumptions on \(\lambda,\mu\) and by the properties of the cut-off function \(\eta\) we obtain
\[
\int_{B_R}(2\lambda|\psi|^2+\big(2-n)\mu|\psi|^4)dx\leq\int_{\mathbb{R}^n}(2\lambda|\psi|^2+\big(2-n)\mu|\psi|^4)\eta(r)d\vol_g
\]
such that we get
\[
\int_{B_R}(2\lambda|\psi|^2+\big(2-n)\mu|\psi|^4)dx\leq C\int_{B_{2R}\setminus B_R}(|\psi||\nabla\psi|+|\psi|^4)dx
\leq C\int_{B_{2R}\setminus B_R}(|\nabla\psi|^\frac{4}{3}+|\psi|^4)dx.
\]
Taking the limit \(R\to\infty\) and making use of the finite energy assumption we obtain
\[
\int_{\mathbb{R}^n}|\psi|^2(2\lambda+(2-n)\mu|\psi|^2)d\vol_g\leq 0,
\]
yielding the result. 
The case \(\lambda\leq 0\) and \(\mu\geq 0\) follows similarly.
By applying the Theorem of Cartan-Hadamard the proof carries over to hyperbolic space.
\end{proof}

\begin{Bem}
In particular, the last Proposition applies in the case \(\mu=0\), 
which corresponds to \(\psi\) being an eigenspinor with eigenvalue \(\lambda\).
Thus, there does not exist an eigenspinor satisfying
\[
\int_M(|\psi|^4+|\nabla\psi|^\frac{4}{3})d\vol_g<\infty
\]
with eigenvalue \(\lambda\) on \(M=\mathbb{R}^n, \mathbb{H}^n\) for \(n\geq 3\).
\end{Bem}

\subsection{Monotonicity formulas for smooth solutions}
In this section we will derive a monotonicity formula for smooth solutions of \eqref{soler} on
complete Riemannian manifolds. We will make use of the fact that the stress-energy tensor \eqref{stress-energy-soler}
is divergence free, whenever \(\psi\) is a solution of \eqref{soler}.
First of all, let us recall the following facts:
A vector field \(X\) is called \emph{conformal} if
\[
\mathcal{L}_Xg=fg,
\]
where \(\mathcal{L}\) denotes the Lie-derivative of the metric with respect to \(X\) and
\(f\colon M\to\mathbb{R}\) is a smooth function.
\begin{Lem}
Let \(T\) be a symmetric 2-tensor. For any vector field \(X\) the following formula holds
\begin{equation}
\operatorname{div}(\iota_X T)=\iota_X\operatorname{div} T+\langle T,\nabla X\rangle.
\end{equation}
If \(X\) is a conformal vector field, then the second term on the right hand side acquires the form
\begin{equation}
\label{conformal-vf}
\langle T,\nabla X\rangle=\frac{1}{n}\operatorname{div}X\tr T.
\end{equation}
\end{Lem}
By integrating over a compact region \(U\), making use of Stokes theorem, we obtain
\begin{Lem}
Let \((M,g)\) be a Riemannian manifold and \(U\subset M\) be a compact region with smooth boundary.
Then, for any symmetric \(2\)-tensor and any vector field \(X\) the following formula holds
\begin{equation}
\label{gauss-tensor-formula}
\int_{\partial U}T(X,\nu)d\sigma=\int_U\iota_X\operatorname{div} Tdx+\int_U\langle T,\nabla X\rangle dx,
\end{equation}
where \(\nu\) denotes the normal to \(U\).
The same formula holds for a conformal vector field \(X\) if we replace the second 
term on the right hand by \eqref{conformal-vf}.
\end{Lem}

We now derive a type of monotonicity formula for smooth solutions of \eqref{soler} in \(\mathbb{R}^n\).

\begin{Prop}[Monotonicity formula in \(\mathbb{R}^n\)]
\label{monotonicity-type-rn}
Let \(\psi\) be a smooth solution of \eqref{soler} on \(M=\mathbb{R}^n\). 
Let \(B_R(x_0)\) be a geodesic ball around the point \(x_0\in M\) and \(0<R_1<R_2\leq R\).
Then the following monotonicity formula holds
\begin{align}
\label{pre-monotonicity}
R_2^{2-n}\mu\int_{B_{R_2(x_0)}}|\psi|^4dx-R_1^{2-n}\mu\int_{B_{R_1(x_0)}}|\psi|^4dx=&-2\lambda\int_{R_1}^{R_2}\big(r^{1-n}\int_{B_r(x_0)}|\psi|^2dx\big)dr \\
&\nonumber+2\int_{R_1}^{R_2}\big(r^{2-n}\int_{\partial B_r(x_0)}\langle\psi,\partial_r\cdot\nabla_{\partial_r}\psi\rangle d\sigma\big) dr.
\end{align}
\end{Prop}

\begin{proof}
For \(M=\mathbb{R}^n\) we choose the conformal vector field \(X=r\frac{\partial}{\partial r}\) with \(r=|x|\).
In this case we have \(\operatorname{div}(X)=n\), thus 
\[
(2-n)\mu\int_{B_r}|\psi|^4dx+r\mu\int_{\partial B_r}|\psi|^4d\sigma=
-2\lambda\int_{B_r}|\psi|^2dx+2r\int_{\partial B_r}\langle\psi,\partial_r\cdot\nabla_{\partial_r}\psi\rangle d\sigma,
\]
where we used \eqref{conformal-vf} and \eqref{gauss-tensor-formula}.
Making use of the coarea formula we can rewrite this as 
\[
\frac{d}{dr}\big(r^{2-n}\mu\int_{B_r}|\psi|^4dx\big)=-2\lambda r^{1-n}\int_{B_r}|\psi|^2dx+2r^{2-n}\int_{\partial B_r}\langle\psi,\partial_r\cdot\nabla_{\partial_r}\psi\rangle d\sigma
\]
and integrating with respect to \(r\) yields the result.
\end{proof}

\begin{Bem}
The previous monotonicity formula also holds if \(\psi\) was only a weak solution of \eqref{soler}, 
that is \(\psi\in L^4(\Sigma M)\times W^{1,\frac{4}{3}}(\Sigma M)\).
\end{Bem}

We now aim at generalizing the monotonicity formula \eqref{pre-monotonicity} to the case of a complete Riemannian
spin manifold. Note that, in general, the vector field \(X=r\frac{\partial}{\partial r}\) will not be conformal.
We fix a point \(x_0\in M\) and consider a ball with geodesic radius \(r=d(x_0,\cdot)\) around that point, 
where \(d\) denotes the Riemannian distance function. Moreover, \(i_M\) will refer to the injectivity radius of \(M\).
Using geodesic polar coordinates we decompose the metric in \(B_{i_M}\) with the help of the Gauss Lemma as
\[
g=g_r+dr\otimes dr.
\]
In the following we will frequently make use of the Hessian of the Riemannian distance function.
Since the Hessian is a symmetric bilinear form we may diagonalize it, its eigenvalues will be denoted by \(\omega_i,i=1,\ldots,n\).
Thus, we may write 
\begin{equation}
\label{Omega}
\sum_{i=1}^n\hess(r^2)(e_i,e_i)=\sum_{i=1}^n\omega_i:=\Omega.
\end{equation}
We denote its largest eigenvalue by \(\omega_{max}\).
The eigenvalues of the Hessian of the Riemannian distance function depend on the geometry of the manifold \(M\) and, in general,
they cannot be computed explicitly.
For some explicit estimates on \(\Omega\) in terms of geometric data we refer to \cite{MR2929724}, Lemma 3.2.

\begin{Lem}
Let \((M,g)\) be a complete Riemannian spin manifold and suppose that \(\psi\) is a smooth solution of \eqref{soler}.
Then the following formula holds
\begin{align}
\label{pre-mono-soler-manifold}
\nonumber
(2\omega_{max}-\Omega)\mu\int_{B_r}|\psi|^4dx+r\mu\int_{\partial B_r}|\psi|^4d\sigma
=&2r\int_{\partial B_r}\langle\psi,\partial_r\cdot\nabla_{\partial_r}\psi\rangle d\sigma 
-2\omega_{max}\lambda\int_{B_r}|\psi|^2dx \\
&-2\sum_{j=2}^n\int_{B_r}\langle e_j\cdot\nabla_{e_j}\psi,\psi\rangle(\omega_j-\omega_{max})dx.
\end{align}
\end{Lem}

\begin{proof}
Inserting the stress-energy tensor \eqref{stress-energy-soler} into \eqref{gauss-tensor-formula} and
choosing the vector field \(X=r\frac{\partial}{\partial r}\) we obtain the following equation
\begin{align*}
2r\int_{\partial B_r}\langle\psi,\partial_r\cdot\nabla_{\partial_r}\psi\rangle d\sigma
-r\mu\int_{\partial B_r}|\psi|^4d\sigma=&\int_{B_r}\langle e_i\cdot\nabla_{e_j}\psi+e_j\cdot\nabla_{e_i}\psi,\psi\rangle\hess(r^2)(e_i,e_j)dx \\
&-\mu\int_{B_r}\tr\hess(r^2)|\psi|^4dx.
\end{align*}
Without loss of generality we assume that \(\omega_1=\omega_{max}\) is the largest eigenvalue of \(\hess(r^2)\).
Diagonalizing the Hessian of the Riemannian distance function we may rewrite
\begin{align*}
\langle e_i\cdot\nabla_{e_j}\psi+e_j\cdot\nabla_{e_i}\psi,\psi\rangle\hess(r^2)(e_i,e_j)=&
2\omega_{max}\langle\psi,D\psi\rangle+2\sum_{j=2}^n\langle e_j\cdot\nabla_{e_j}\psi,\psi\rangle(\omega_j-\omega_{max})\\
=&2\omega_{max}(\lambda|\psi|^2+\mu|\psi|^4)\\
&+2\sum_{j=2}^n\langle e_j\cdot\nabla_{e_j}\psi,\psi\rangle(\omega_j-\omega_{max}),
\end{align*}
which yields the claim.
\end{proof}

\begin{Bem}
The problematic contributions in the monotonicity-type formulas \eqref{monotonicity-type-rn} and \eqref{pre-mono-soler-manifold}
are the indefinite terms \(\langle\psi,\partial_r\cdot\nabla_{\partial_r}\psi\rangle\) and 
\(\langle e_i\cdot\nabla_{e_i}\psi,\psi\rangle\hess(r^2)(e_i,e_i)\).
To give them a definite sign we could assume that \(\psi\) is both a solution of \eqref{soler} and a twistor spinor. 
In this case we would have
\[
\langle\psi,\partial_r\cdot\nabla_{\partial_r}\psi\rangle=\frac{1}{n}g(\partial_r,\partial_r)\langle\psi,D\psi\rangle
=\frac{1}{n}g(\partial_r,\partial_r)(\lambda|\psi|^2+\mu|\psi|^4).
\]
The right hand side of this equation is positive for \(\lambda,\mu>0\).
However, we have already seen that under the assumptions from above \(|\psi|^2\) 
is equal to a constant and in this case the monotonicity formula contains no interesting information.
Moreover, regarding the second term, we would get
\begin{align*}
\langle e_i\cdot\nabla_{e_i}\psi,\psi\rangle\hess(r^2)(e_i,e_i)=&-\frac{1}{n}\langle e_i\cdot D\psi,\psi\rangle\hess(r^2)(e_i,e_i) \\
=&-\frac{1}{n}\langle e_i\cdot\psi,\psi\rangle(\lambda|\psi|^2+\mu|\psi|^4)\hess(r^2)(e_i,e_i)=0.
\end{align*}
\end{Bem}

\begin{Bem}
It would be desirable to estimate the term \(\langle e_i\cdot\nabla_{e_i}\psi,\psi\rangle\hess(r^2)(e_i,e_i)\) in \eqref{pre-mono-soler-manifold}
in terms of geometric data of the manifold \(M\) and the right hand side of \eqref{soler}. Unfortunately, this only seems to be possible
if all eigenvalues of the Hessian of the Riemann distance function would be equal.
\end{Bem}

\begin{Prop}
Let \((M,g)\) be a complete Riemannian spin manifold and suppose that \(\psi\) is a smooth solution of \eqref{soler}.
Then for all \(0<R_1<R_2\leq R\), \(R\in(0,i_M)\) the following type of monotonicity formula holds
\begin{align}
\label{mono-soler-manifold}
R_1^{2\omega_{max}-\Omega}&\int_{B_{R_1}}(\mu|\psi|^4-2\langle\psi,\partial_r\cdot\nabla_{\partial_r}\psi\rangle) dx \\
\nonumber =&R_2^{2\omega_{max}-\Omega}\int_{B_{R_2}}(\mu|\psi|^4-2\langle\psi,\partial_r\cdot\nabla_{\partial_r}\psi\rangle) dx \\
\nonumber &+2(2\omega_{max}-\Omega)\int_{R_1}^{R_2}\big(r^{{2\omega_{max}-\Omega}-1}\int_{B_{r}}\langle\psi,\partial_r\cdot\nabla_{\partial_r}\psi\rangle dx\big)dr \\
\nonumber &+2\omega_{max}\lambda \int_{R_1}^{R_2}\big(r^{2\omega_{max}-\Omega-1}\int_{B_r}|\psi|^2dx\big)dr \\
\nonumber &+2\sum_{j=2}^n(\omega_j-\omega_{max})\int_{R_1}^{R_2}\big(r^{2\omega_{max}-\Omega-1}\int_{B_r}\langle e_j\cdot\nabla_{e_j}\psi,\psi\rangle\big)dr,
\end{align}
where \(\Omega\) is given by \eqref{Omega}.
\end{Prop}
\begin{proof}
Using \eqref{pre-mono-soler-manifold} and the coarea formula we find
\begin{align*}
\frac{d}{dr}r^{2\omega_{max}-\Omega}\mu\int_{B_r}|\psi|^4dx=&
2r^{2\omega_{max}-\Omega}\int_{\partial B_r}\langle\psi,\partial_r\cdot\nabla_{\partial_r}\psi\rangle d\sigma
-2\omega_{max}\lambda r^{2\omega_{max}-\Omega-1}\int_{B_r}|\psi|^2dx \\
&-2\sum_{j=2}^nr^{2\omega_{max}-\Omega-1}\int_{B_r}\langle e_j\cdot\nabla_{e_j}\psi,\psi\rangle(\omega_j-\omega_{max})dx.
\end{align*}
Integrating with respect to \(r\) and using integration by parts 
\begin{align*}
\int_{R_1}^{R_2}&\big(r^{2\omega_{max}-\Omega}\int_{\partial B_r}\langle\psi,\partial_r\cdot\nabla_{\partial_r}\psi\rangle d\sigma\big)dr \\
=&\int_{R_1}^{R_2}\big(r^{2\omega_{max}-\Omega}\frac{d}{dr}\int_{B_r}\langle\psi,\partial_r\cdot\nabla_{\partial_r}\psi\rangle dx\big)dr \\
=&R_2^{2\omega_{max}-\Omega}\int_{B_{R_2}}\langle\psi,\partial_r\cdot\nabla_{\partial_r}\psi\rangle dx 
-R_1^{2\omega_{max}-\Omega}\int_{B_{R_1}}\langle\psi,\partial_r\cdot\nabla_{\partial_r}\psi\rangle dx \\
&+(\Omega-2\omega_{max})\int_{R_1}^{R_2}\big(r^{{2\omega_{max}-\Omega}-1}\int_{B_{r}}\langle\psi,\partial_r\cdot\nabla_{\partial_r}\psi\rangle dx\big)dr
\end{align*}
yields the claim.
\end{proof}

\begin{Bem}
If \(M=\mathbb{R}^n\), then \(\omega_i=1,i=1,\ldots,n\) and \(\Omega=n\).
In this case \eqref{mono-soler-manifold} reduces to \eqref{monotonicity-type-rn}.
\end{Bem}

\begin{Bem}
It seems very difficult to obtain a Liouville Theorem from the monotonicity formula \eqref{mono-soler-manifold}
without posing many conditions on the solution of \eqref{soler}.
\end{Bem}

\subsection{A Liouville Theorem for complete manifolds with positive Ricci curvature}
In this section we will prove a Liouville theorem for smooth solutions of \eqref{soler} on complete 
noncompact manifolds with positive Ricci curvature. Our result is motivated
from a similar result for harmonic maps, see \cite{MR0438388}, Theorem 1.
We set \(e(\psi):=\frac{1}{2}|\psi|^4\).
\begin{Satz}
\label{theorem-liouville-complete-soler}
Let \((M,g)\) be a complete noncompact Riemannian spin manifold with positive Ricci curvature.
Suppose that 
\begin{equation}
\label{assumption-liouville-complete-soler}
R\geq 4(\lambda+\mu|\psi|^2)^2.
\end{equation}
If \(\psi\) is a smooth solution of \eqref{soler} with finite energy \(e(\psi)\) then \(\psi\) vanishes identically.
\end{Satz}

\begin{proof}
Making use of the assumption the Bochner formula \eqref{bochner} yields
\begin{equation}
\label{liouville-complete-a}
\Delta e(\psi)\geq \big|d|\psi|^2\big|^2.
\end{equation}
In addition, by the Cauchy-Schwarz inequality we find
\begin{equation}
\label{liouville-complete-b}
|de(\psi)|^2\leq 2e(\psi)\big|d|\psi|^2\big|^2.
\end{equation}
We fix a positive number \(\epsilon>0\) and calculate
\begin{align*}
\Delta\sqrt{e(\psi)+\epsilon}=\frac{\Delta e(\psi)}{2\sqrt{e(\psi)+\epsilon}}
-\frac{1}{4}\frac{|de(\psi)|^2}{(e(\psi)+\epsilon)^\frac{3}{2}} 
\geq\frac{\big|d|\psi|^2\big|^2}{2\sqrt{e(\psi)+\epsilon}}\big(1-\frac{e(\psi)}{e(\psi)+\epsilon}\big)
\geq 0,
\end{align*}
where we used \eqref{liouville-complete-a} and \eqref{liouville-complete-b}.
Let \(\eta\) be an arbitrary function on \(M\) with compact support. 
We obtain
\begin{align*}
0\leq&\int_M\eta^2\sqrt{e(\psi)+\epsilon}\Delta\sqrt{e(\psi)+\epsilon}d\vol_g \\
=&-2\int_M\eta\sqrt{e(\psi)+\epsilon}\langle\nabla\eta,\nabla\sqrt{e(\psi)+\epsilon}\rangle d\vol_g
-\int_M\eta^2|\nabla\sqrt{e(\psi)+\epsilon}|^2d\vol_g.
\end{align*}
Now let \(x_0\) be a point in \(M\) and let \(B_R,B_{2R}\) be geodesic balls centered at \(x_0\)
with radii \(R\) and \(2R\).
We choose a cutoff function \(\eta\) satisfying
\[
\eta(x)=
\begin{cases}
1,\qquad x\in B_R,\\
0, \qquad x\in M\setminus B_{2R}.
\end{cases}
\]
In addition, we choose \(\eta\) such that
\[
0\leq\eta\leq 1,\qquad |\nabla\eta|\leq\frac{C}{R}
\]
for a positive constant \(C\).
Then, we find
\begin{align*}
0\leq&-2\int_{B_{2R}}\eta\sqrt{e(\psi)+\epsilon}\langle\nabla\eta,\nabla\sqrt{e(\psi)+\epsilon}\rangle dx
-\int_{B_{2R}}\eta^2|\nabla\sqrt{e(\psi)+\epsilon}|^2 dx\\
\leq&2\big(\int_{B_{2R\setminus B_R}}\eta^2|\sqrt{e(\psi)+\epsilon}|^2dx\big)^\frac{1}{2}
\big(\int_{B_{2R\setminus B_R}}|\nabla\eta|^2(e(\psi)+\epsilon)dx\big)^\frac{1}{2}\\
&-\int_{B_{2R}{\setminus B_R}}\eta^2|\nabla\sqrt{e(\psi)+\epsilon}|^2dx
-\int_{B_R}|\nabla\sqrt{e(\psi)+\epsilon}|^2dx.
\end{align*}
We therefore obtain
\begin{align*}
\int_{B_r}|\nabla\sqrt{e(\psi)+\epsilon}|^2dx\leq\int_{B_{2R}{\setminus B_R}}|\nabla\eta|^2(e(\psi)+\epsilon)dx
\leq\frac{C^2}{R^2}\int_{B_{2R}}(e(\psi)+\epsilon)dx.
\end{align*}
We set \(B'_R:=B_R\setminus\{x\in B_R\mid e(\psi)(x)=0\}\) and find
\begin{align*}
\int_{B'_r}\frac{|\nabla(e(\psi)+\epsilon)|^2}{4(e(\psi)+\epsilon)}dx
\leq\frac{C^2}{R^2}\int_{B_{2R}}(e(\psi)+\epsilon)dx.
\end{align*}
Letting \(\epsilon\to 0\) we get
\begin{align*}
\int_{B'_r}\frac{|\nabla(e(\psi)|^2}{4e(\psi)}dx
\leq\frac{C^2}{R^2}\int_{B_{2R}}e(\psi)dx.
\end{align*}
Now, letting \(R\to\infty\) and under the assumption that the energy is finite, we have
\[
\int_{M\setminus\{e(\psi)=0\}}\frac{|\nabla e(\psi)|^2}{4e(\psi)}d\vol_g\leq 0,
\]
hence the energy \(e(\psi)\) has to be constant.
If \(e(\psi)\neq 0\), then the volume of \(M\) would have to be finite.
However, by Theorem 7 of \cite{MR0417452} the volume of a complete and noncompact Riemannian 
manifold with nonnegative Ricci curvature is infinite. Hence \(e(\psi)=0\),
which yields the result.
\end{proof}

Note, that Theorem \ref{theorem-liouville-complete-soler} also holds in the case \(\mu=0\),
which gives us the following vanishing result for eigenspinors:
\begin{Cor}
Suppose that \(\psi\) is a smooth solution of \(D\psi=\lambda\psi\) on 
a complete noncompact manifold with positive Ricci curvature. If
\begin{equation*}
R\geq 4\lambda^2
\end{equation*}
and \(e(\psi)\) is finite then \(\psi\) vanishes identically.
\end{Cor}

\section{Dirac-harmonic maps with curvature term from complete manifolds}
\emph{Dirac-harmonic maps with curvature term} arise as critical points of part of the supersymmetric nonlinear \(\sigma\)-model from
quantum field theory \cite{MR1701618}, p. 268, the only difference being that in contrast to the physics literature standard, that is commuting, spinors are used.
They form a pair of a map from a Riemann spin manifold to another Riemannian manifold coupled with a vector spinor.
For a two-dimensional domain they belong to the class of conformally invariant variational problems.
The conformal invariance gives rise to a removable singularity theorem \cite{MR3558358} and an energy identity \cite{72486}.
Conservation laws for Dirac-harmonic maps with curvature term were established in \cite{MR3735550}
and a vanishing result for the latter under small-energy assumptions was derived in \cite{MR3886921}.
For Dirac-wave maps with curvature term (which are Dirac-harmonic maps with curvature term from a domain with Lorentzian metric)
on expanding spacetimes an existence result could be achieved in \cite{MR3830277}.

The mathematical study of the supersymmetric nonlinear \(\sigma\)-model with standard spinors was initiated in \cite{MR2262709}, 
where the notion of \emph{Dirac-harmonic maps} was introduced.
The full action of the supersymmetric nonlinear \(\sigma\)-model contains two additional terms:
Taking into account an additional two-form in the action functional the resulting equations were studied in \cite{MR3305429},
Dirac-harmonic maps with curvature term to target spaces with torsion are analyzed in \cite{MR3493217}.

Most of the results presented in this section still hold true if we would consider the full supersymmetric nonlinear \(\sigma\)-model.
Let us give some more details in support of this statement: The central ingredient in the derivation of various monotonicity formulas
and Liouville theorems will be the stress-energy tensor. An additional two-form contribution in the action functional
would not give a contribution to the stress-energy tensor as it does not depend on the metric of the domain, see
\cite[Section 3]{MR3305429} for more details. 
Moreover, if we would consider a connection with torsion on the target manifold
we would get the same stress-energy tensor, see \cite[Section 4]{MR3493217}, 
and all results that will be formulated in this section still hold if we 
formulate the curvature assumptions taking into account the connection with torsion.

Let us again emphasize that in the physics literature anticommuting spinors are employed while the mathematical references 
stated above and the present article consider standard commuting spinors.

In the following we still assume that \((M,g)\) is a complete Riemannian spin manifold
and \((N,h)\) another Riemannian manifold. Whenever we will make use of indices
we use Latin letters for indices related to \(M\) and Greek letters for indices related to \(N\).
Let \(\phi\colon M\to N\) be a map and let \(\phi^\ast TN\) be the pull-back of the tangent bundle from \(N\).
We consider the twisted bundle \(\Sigma M\otimes\phi^\ast TN\), on this bundle
we obtain a connection induced from \(\Sigma M\) and \(\phi^\ast TN\), which will be
denoted by \(\tilde{\nabla}\). Sections in \(\Sigma M\otimes\phi^\ast TN\) are called \emph{vector spinors}.
On \(\Sigma M\otimes\phi^\ast TN\) we have a scalar product induced from \(\Sigma M\) and \(\phi^\ast TN\),
we will denote its real part by \(\langle\cdot,\cdot\rangle\).
The twisted Dirac operator acting on vector spinors is defined as
\[
\D:=\sum_{i=1}^ne_i\cdot\tilde{\nabla}_{e_i}.
\]
Note that the operator \(\D\) is still elliptic. Moreover, we assume that the connection on \(\phi^\ast TN\) is metric, thus \(\D\)
is also self-adjoint with respect to the \(L^2\)-norm if \(M\) is compact.
The action functional for Dirac-harmonic maps with curvature term is given by
\begin{equation}
\label{energy-dhc}
E_c(\phi,\psi)=\frac{1}{2}\int_M(|d\phi|^2+\langle\psi,\D\psi\rangle-\frac{1}{6}\langle R^N(\psi,\psi)\psi,\psi\rangle)d\vol_g.
\end{equation}
Here, \(R^N\) denotes the curvature tensor of the manifold \(N\).
The factor \(1/6\) in front of the curvature term is required by supersymmetry, see \cite{MR1701618}.
The indices are contracted as 
\[
\langle R^N(\psi,\psi)\psi,\psi\rangle=R_{\alpha\beta\gamma\delta}\langle\psi^\alpha,\psi^\gamma\rangle\langle\psi^\beta,\psi^\delta\rangle,
\]
which ensures that the functional is real valued.
The critical points of the action functional \eqref{energy-dhc} are given by
\begin{align}
\label{euler-lagrange-phi}\tau(\phi)=&\frac{1}{2}R^N(\psi,e_i\cdot\psi)d\phi(e_i)
-\frac{1}{12}\langle(\nabla R^N)^\sharp(\psi,\psi)\psi,\psi\rangle, \\
\label{euler-lagrange-psi}\D\psi=&\frac{1}{3}R^N(\psi,\psi)\psi,
\end{align}
where \(\tau(\phi)\in\Gamma(\phi^\ast TN)\) is the tension field of the map \(\phi\) and \(\sharp\colon\phi^\ast T^\ast N\to\phi^\ast TN\) represents the musical isomorphism.
For a derivation see \cite{MR2370260}, Section II and \cite{MR3333092}, Proposition 2.1.

Solutions \((\phi,\psi)\) of the system \eqref{euler-lagrange-phi}, \eqref{euler-lagrange-psi}
are called \emph{Dirac-harmonic maps with curvature term}.

The correct function space for weak solutions of \eqref{euler-lagrange-phi}, \eqref{euler-lagrange-psi} is 
\[
\chi(M,N):=W^{1,2}(M,N)\times W^{1,\frac{4}{3}}(M,\Sigma M\otimes\phi^\ast TN)\times L^4(M,\Sigma M\otimes\phi^\ast TN).
\]
For the domain being a closed surface it was shown in \cite{MR3333092} that a weak solution \((\phi,\psi)\in\chi(M,N)\) of \eqref{euler-lagrange-phi}, \eqref{euler-lagrange-psi}
is smooth. This was later extended to higher dimensions in \cite{MR4018319}, see also \cite{MR3798022} for the regularity
of Dirac-harmonic maps with curvature term coupled to a gravitino.

For smooth solutions of \eqref{euler-lagrange-phi}, \eqref{euler-lagrange-psi} on a closed Riemannian surface
a vanishing result was obtained in \cite{MR3333092}, Lemma 4.9.
More precisely, it was shown that a smooth Dirac-harmonic map with curvature term with small energy \(\int_M(|d\phi|^2+|\psi|^4)d\vol_g\)
from a closed surface that does not admit ``standard'' harmonic spinors must be trivial.
Using the recent regularity for vector spinors \cite{MR3908762}, Theorem 1.2 it should be possible
to prove this result without the assumption that \(M\) is not allowed to have harmonic spinors.
However, Theorem 1.2 in \cite{MR3908762} is formulated for the case of a domain manifold with boundary
and one would require a version for closed manifolds.

\begin{Dfn}
A weak Dirac-harmonic map with curvature term \((\phi,\psi)\in\chi(M,N)\) 
is called \emph{stationary} if it is also a critical point of \(E_c(\phi,\psi)\) with respect
to domain variations.
\end{Dfn}

To obtain the formula for stationary Dirac-harmonic maps with curvature term we make use of 
the same methods as before. Since the twist bundle \(\phi^\ast TN\) 
does not depend on the metric on \(M\) we can use the same methods as in Section 2.
Thus, let \(k\) be a smooth element of Sym\((0,2)\).
Again, we will use the notation \(\psi_{g+tk}:=\beta_{g,g+tk}\psi\in\Gamma(\Sigma M_{g+tk}\otimes\phi^\ast TN)\).

\begin{Lem}
The following formula for the variation of the twisted Dirac-energy with respect to the metric holds
\begin{align}
\label{variation-dirac-dhc}
\frac{d}{dt}\big|_{t=0}\langle\psi_{g+tk},\D_{g+tk}&\psi_{g+tk}\rangle_{\Sigma_{g+tk}M\otimes\phi^\ast TN} \\
\nonumber=&-\frac{1}{4}\langle e_i\cdot\nabla^{\Sigma_gM\otimes\phi^\ast TN}_{e_j}\psi
+e_j\cdot\nabla^{\Sigma_gM\otimes\phi^\ast TN}_{e_i}\psi,\psi\rangle_{\Sigma_gM\otimes\phi^\ast TN}k^{ij}
\end{align}
with the stress-energy tensor associated to the twisted Dirac energy on the right hand side.
\end{Lem}

At this point we are ready to compute the variation of the action functional \eqref{energy-dhc} with respect to the metric.

\begin{Prop}
Let the pair \((\phi,\psi)\in \chi(M,N)\) be a weak Dirac-harmonic map with curvature term.
Then \((\phi,\psi)\) is a stationary Dirac-harmonic map with curvature term if for any smooth symmetric \((0,2)\)-tensor \(k\)
the following formula holds
\begin{align}
\label{stationary-identity-dhc}
\int_M \big(&2\langle d\phi(e_i),d\phi(e_j)\rangle-g_{ij}|d\phi|^2
+\frac{1}{2}\langle\psi,e_i\cdot\nabla^{\Sigma_gM\otimes\phi^\ast TN}_{e_j}\psi+e_j\cdot\nabla^{\Sigma_gM\otimes\phi^\ast TN}_{e_i}\psi\rangle \\
\nonumber&-\frac{1}{6}g_{ij}\langle R^N(\psi,\psi)\psi,\psi\rangle\big)k^{ij})d\vol_g=0.
\end{align}
\end{Prop}

\begin{proof}
We calculate
\[
\frac{d}{dt}\big|_{t=0} E_c(\phi,\psi,g+tk)=0,
\]
where \(k\) is a symmetric \((0,2)\)-tensor and \(t\) some small number.
Using the variation of the volume-element \eqref{variation-volume} we obtain the variation of the Dirichlet energy 
\[
\frac{d}{dt}\big|_{t=0}\int_M|d\phi|^2_{g+tk}d\vol_{g+tk}
=\int_M\big(-\langle h(d\phi(e_i),d\phi(e_j)),k_{ij}\rangle d\vol_g
+\frac{1}{2}|d\phi|^2\langle g,k\rangle_gd\vol_g\big).
\]
Note that we get a minus sign in the first term since \(d\phi\in\Gamma(T^\ast M\otimes\phi^\ast TN)\)
such that we have to vary the metric on the cotangent bundle.
As a second step, we compute the variation of the Dirac energy using \eqref{variation-dirac-dhc} and \eqref{variation-volume} yielding
\begin{align*}
\frac{d}{dt}&\big|_{t=0}\int_M\langle\psi_{g+tk},\D_{g+tk}\psi_{g+tk}\rangle d\vol_{g+tk}\\
=&\int_M-\frac{1}{4}\langle e_i\cdot\nabla^{\Sigma_gM\otimes\phi^\ast TN}_{e_j}\psi
+e_j\cdot\nabla^{\Sigma_gM\otimes\phi^\ast TN}_{e_i}\psi,\psi\rangle_{\Sigma_gM}k^{ij}d\vol_g
+\frac{1}{2}\int_M\langle\psi,\D\psi\rangle\langle g,k\rangle_gd\vol_g.
\end{align*}
Finally, for the term involving the curvature tensor of the target and the four spinors we obtain
\begin{align*}
\frac{d}{dt}\big|_{t=0}\int_M\langle R^N(\psi_{g+tk},\psi_{g+tk})\psi_{g+tk},&\psi_{g+tk}\rangle_{g+tk}d\vol_{g+tk} \\
=&\frac{d}{dt}\big|_{t=0}\int_M\langle R^N(\psi,\psi)\psi,\psi\rangle_{\Sigma_{g+tk}M\otimes\phi^\ast TN} d\vol_{g+tk} \\
=&\frac{1}{2}\int_M\langle R^N(\psi,\psi)\psi,\psi\rangle_{\Sigma_gM\otimes\phi^\ast TN}\langle g,k\rangle_gd\vol_g,
\end{align*}
where we used that \(\beta\) acts as an isometry on the spinor bundle in the first step.
Adding up the three contributions and using the fact that \((\phi,\psi)\) is a weak Dirac-harmonic map
with curvature term yields the result.
\end{proof}

\subsection{A Liouville Theorem for stationary solutions}
It is well known that a stationary harmonic map \(\mathbb{R}^q\to N\) with finite Dirichlet energy 
is a constant map \cite{MR539148}, Section 5. This result was generalized to stationary Dirac-harmonic maps and here
we generalize it to stationary Dirac-harmonic maps with curvature term by adding a curvature assumption.

A similar result for smooth Dirac-harmonic maps with curvature term was already obtained in \cite{MR2370260}, Theorem 1.2.
Let us point out in some more detail the similarities and differences between the methods of proof used in \cite{MR2370260}
and in the present article. In the proof of Theorem 1.2 in \cite{MR2370260} the authors calculate the Lie-derivative of
the energy density of \eqref{energy-dhc} with respect to a conformal vector field \(X\).
In order to carry out the Lie-derivative of the terms involving spinors in \eqref{energy-dhc} they also
apply the methods of \cite{MR1158762}. After having obtained a formula for the Lie-derivative of
the energy density of \eqref{energy-dhc} they multiply it with a suitable cutoff function
and the result follows after integration by parts. Although our method of proof formally 
looks very different it has the same core ideas. At its heart is
on the one hand the stress-energy tensor which was also derived using the methods of \cite{MR1158762}
and on the other hand we also crucially require the existence of a conformal vector field.
However, it seems that the advantage of our method is that we do not require to have a smooth solution of \eqref{euler-lagrange-phi}, \eqref{euler-lagrange-psi}.
On the other hand both proofs require the existence of a conformal vector field such that they can only
work on Riemannian manifolds with a sufficient amount of symmetry.

First, we will give the following remark following the proof of Theorem 3.1 in \cite{MR3886921}.
\begin{Bem}
\label{rem-monotone-quantity-dhc}
In this section we will often consider the quantity
\begin{align}
\label{monotone-quantity-dhc}
|d\phi|^2+\frac{1}{6}\langle R^N(\psi,\psi)\psi,\psi\rangle
\end{align}
and it will be crucial for our arguments that this expression is positive.
\begin{enumerate}
 \item In the case that \(\phi\colon M\to N\) is a constant map we can consider \(v\in\phi^\ast TN,\Psi\in\Gamma(\Sigma M)\) and set \(\psi:=\Psi\otimes v\).
It is easy to check that this pair \((\phi,\psi)\) satisfies
\begin{align*}
\langle R^N(\psi,\psi)\psi,\psi\rangle=\langle R^N(v,v)v,v\rangle|\Psi|^4=0
\end{align*}
due to the skew symmetry of the Riemann curvature tensor regardless of any curvature assumptions on the target.
Hence, in this case the system \eqref{euler-lagrange-phi}, \eqref{euler-lagrange-psi} would reduce to
\begin{align*}
D\psi^\alpha=0,\qquad 1\leq\alpha\leq\dim N,
\end{align*}
where \(D\) denotes the standard Dirac operator on \(\Sigma M\).

\item However, for a pair \((\phi,\psi)\) that is not of the form from above 
the term \(\langle R^N(\psi,\psi)\psi,\psi\rangle\) will be different from zero.
A careful inspection reveals that for \(N\) having positive sectional curvature we have 
\begin{align*}
|d\phi|^2+\frac{1}{6}\langle R^N(\psi,\psi)\psi,\psi\rangle\geq 0,
\end{align*}
see \cite[Proof of Theorem 1.2]{MR2370260} for more details.
\end{enumerate}
\end{Bem}

\begin{Satz}
\label{theorem-dhc-stationary}
Let \(M=\mathbb{R}^n,\mathbb{H}^n, n\geq 3\) and suppose that 
\((\phi,\psi)\in W^{1,2}_{loc}(M,N)\times W^{1,\frac{4}{3}}_{loc}(M,\Sigma M\otimes\phi^\ast TN)\times L^{4}_{loc}(M,\Sigma M\otimes\phi^\ast TN)\) 
is a stationary Dirac-harmonic maps with curvature term satisfying
\begin{equation}
\int_{\mathbb{R}^n}(|d\phi|^2+|\nabla^{\Sigma M}\psi|^\frac{4}{3}+|\psi|^4)d\vol_g<\infty. 
\end{equation}
If \(N\) has positive sectional curvature then \(\phi\) is constant and \(\psi\) vanishes identically.
\end{Satz}

\begin{proof}
Let \(\eta\in C_0^\infty(\mathbb{R})\) be a smooth cut-off function satisfying \(\eta=1\) for \(r\leq R\),
\(\eta=0\) for \(r\geq 2R\) and \(|\eta'(r)|\leq\frac{C}{R}\). In addition, 
we choose \(Y(x):=x\eta(r)\in C_0^\infty(\mathbb{R}^n,\mathbb{R}^n)\) with \(r=|x|\).
Hence, we find
\[
k_{ij}=\frac{\partial Y_i}{\partial x^j}=\delta_{ij}\eta(r)+\frac{x_i x_j}{r}\eta'(r).
\]
Inserting this into \eqref{stationary-identity-dhc} and using that \((\phi,\psi)\) is a weak solution of the system \eqref{euler-lagrange-phi}, \eqref{euler-lagrange-psi}
we obtain
\begin{align*}
(2-n)\int_{\mathbb{R}^n}(|d\phi|^2+\frac{1}{6}\langle R^N(\psi,\psi)\psi,\psi\rangle)\eta(r)d\vol_g
=&\int_{\mathbb{R}^n}(|d\phi|^2-2\big|\frac{\partial\phi}{\partial r}\big|^2-\langle\psi,\partial_r\cdot\tilde{\nabla}_{\partial_r}\psi\rangle \\
&+\frac{1}{6}\langle R^N(\psi,\psi)\psi,\psi\rangle)r\eta'(r)d\vol_g.
\end{align*}
By the properties of the cut-off function \(\eta\) we find (see the proof of Theorem \ref{theorem-liouville-stationary-soler} for more details)
\begin{align*}
(2-n)\int_{\mathbb{R}^n}(|d\phi|^2+\frac{1}{6}\langle R^N(\psi,\psi)\psi,\psi\rangle)\eta(r)d\vol_g
\leq & C\int_{B_{2R}\setminus B_R}(|d\phi|^2+|\psi||\tilde{\nabla}\psi|+|\psi|^4)dx \\
\leq & C\int_{B_{2R}\setminus B_R}(|d\phi|^2+|\nabla^{\Sigma M}\psi|^\frac{4}{3}+|\psi|^4)dx.
\end{align*}
Due to the finite energy assumption and the fact that \(n\geq 3\), taking the limit \(R\to\infty\) yields
\[
\int_{\mathbb{R}^n}(|d\phi|^2+\frac{1}{6}\langle R^N(\psi,\psi)\psi,\psi\rangle)d\vol_g=0.
\]
At this point we need to make a case distinction as in Remark \ref{rem-monotone-quantity-dhc}.
In the first case the statement follows from Theorem \ref{theorem-liouville-stationary-soler} with \(\lambda=\mu=0\)
and in the second case we are done since \(N\) has positive sectional curvature.
To obtain the result for hyperbolic space we again apply the theorem of Cartan-Hadamard.
\end{proof}

\subsection{Monotonicity formulas and Liouville Theorems}
In this section we derive a monotonicity formula for Dirac-harmonic maps with curvature term
building on their stress-energy tensor. For simplicity, we will mostly assume that \((\phi,\psi)\) is a smooth Dirac-harmonic map with curvature term.
From \eqref{stationary-identity-dhc} we obtain the stress-energy tensor for the functional \(E_c(\phi,\psi)\) as
\begin{align}
\label{stress-energy-dhc}
S_{ij}=&2\langle d\phi(e_i),d\phi(e_j)\rangle-g_{ij}|d\phi|^2 \\
&+\frac{1}{2}\langle\psi,e_i\cdot\nabla^{\Sigma M\otimes\phi^\ast TN}_{e_j}\psi+e_j\cdot\nabla^{\Sigma M\otimes\phi^\ast TN}_{e_i}\psi\rangle 
\nonumber-\frac{1}{6}g_{ij}\langle R^N(\psi,\psi)\psi,\psi\rangle.
\end{align}
It is well-known that the stress-energy tensor \eqref{stress-energy-dhc} is divergence free
in the case of a two-dimensional domain whenever \((\phi,\psi)\) solves the equation for Dirac-harmonic maps with
curvature term.
This question was first addressed in \cite{MR3333092}, Proposition 3.2. 
However, in the calculation carried out in that reference a real-part in front of the third term is missing.
This issue was later clarified and corrected in \cite{72486}, Lemma 4.1.

For the sake of completeness and in order to also include the case of a higher-dimensional domain manifold
we will give another proof that \eqref{stress-energy-dhc} is divergence free.

\begin{Lem}
Suppose that \((\phi,\psi)\) is a smooth solution of \eqref{euler-lagrange-phi}, \eqref{euler-lagrange-psi}.
Then the stress-energy tensor \eqref{stress-energy-dhc} is divergence-free.
\end{Lem}
\begin{proof}
First, we replace the last term in \eqref{stress-energy-dhc} using \eqref{euler-lagrange-psi}.
To shorten the notation we will write \(\tilde\nabla\) for the connection on \(\Sigma M\otimes\phi^\ast TN\).
Then the stress-energy tensor acquires the form
\begin{align*}
S_{ij}=2\langle d\phi(e_i),d\phi(e_j)\rangle-g_{ij}|d\phi|^2 
+\frac{1}{2}\langle\psi,e_i\cdot\tilde\nabla_{e_j}\psi+e_j\cdot\tilde\nabla_{e_i}\psi\rangle 
-\frac{1}{2}g_{ij}\langle\D\psi,\psi\rangle.
\end{align*}
We choose a local orthonormal basis of \(TM\) such that \(\nabla_{e_i}e_j=0,i,j=1,\ldots,n\) at the considered point.
By a direct calculation we find
\begin{align}
\label{identity-a}
\nabla^j\big(2\langle d\phi(e_i),d\phi(e_j)\rangle-g_{ij}|d\phi|^2\big)=&2\langle d\phi(e_i),\tau(\phi)\rangle \\
\nonumber=&\langle d\phi(e_i),R^N(\psi,e_r\cdot\psi)d\phi(e_r)\rangle \\
\nonumber&-\frac{1}{6}\big\langle d\phi(e_i),\langle(\nabla R^N)^\sharp(\psi,\psi)\psi,\psi\rangle\big\rangle,
\end{align}
where we have used \eqref{euler-lagrange-phi} in the second step.
Then, we calculate
\begin{align}
\label{identity-b}
\nabla^j\big(\frac{1}{2}&\langle\psi,e_i\cdot\tilde\nabla_{e_j}\psi+e_j\cdot\tilde\nabla_{e_i}\psi\rangle 
-\frac{1}{2}g_{ij}\langle \D\psi,\psi\rangle\big) \\
\nonumber=&\frac{1}{2}\underbrace{\langle\tilde\nabla_{e_j}\psi,e_i\cdot\tilde\nabla_{e_j}\psi\rangle}_{=0}
+\frac{1}{2}\langle\psi,e_i\cdot\tilde\Delta\psi\rangle 
+\frac{1}{2}\langle\tilde\nabla_{e_j}\psi,e_j\cdot\tilde\nabla_{e_i}\psi\rangle
+\frac{1}{2}\langle\psi,\D\tilde\nabla_{e_i}\psi\rangle \\
\nonumber&-\frac{1}{2}\langle\tilde\nabla_{e_i}\psi,\D\psi\rangle
-\frac{1}{2}\langle\psi,\tilde\nabla_{e_i}\D\psi\rangle \\
\nonumber=&\frac{1}{2}\langle\psi,e_i\cdot\tilde\Delta\psi\rangle
-\langle\D\psi,\tilde\nabla_{e_i}\psi\rangle 
+\frac{1}{2}\langle\psi,\D\tilde\nabla_{e_i}\psi-\tilde\nabla_{e_i}\D\psi\rangle.
\end{align}
Recall that
\begin{align*}
\D\tilde\nabla_{e_i}\psi=e_r\cdot R^{\Sigma M}(e_r,e_i)\psi+e_r\cdot R^N(d\phi(e_r),d\phi(e_i))\psi+\tilde\nabla_{e_i}\D\psi
\end{align*}
such that 
\begin{align}
\label{identity-c}
\langle\psi,\D\tilde\nabla_{e_i}\psi-\tilde\nabla_{e_i}\D\psi\rangle=\frac{1}{2}\underbrace{\langle\psi,\operatorname{Ric}(e_i)\cdot\psi\rangle}_{=0}
+\langle\psi,e_r\cdot R^N(d\phi(e_r),d\phi(e_i))\psi\rangle.
\end{align}
In order to manipulate the term involving the connection Laplacian on \(\Sigma M\otimes\phi^\ast TN\) we recall
the Weitzenböck formula for the twisted Dirac operator \(\D\) which is given by
\begin{equation*}
\D^2\psi=-\tilde{\Delta}\psi+\frac{1}{4}R\psi+\frac{1}{2}e_r\cdot e_s\cdot R^N(d\phi(e_r),d\phi(e_s))\psi.
\end{equation*}
This allows us to conclude that
\begin{align*}
\langle\psi,e_i\cdot\tilde\Delta\psi\rangle=-\langle\psi,e_i\cdot\D^2\psi\rangle+\frac{R}{4}\underbrace{\langle\psi,e_i\cdot\psi\rangle}_{=0}
+\frac{1}{2}\langle\psi,e_i\cdot e_r\cdot e_s\cdot R^N(d\phi(e_r),d\phi(e_s))\psi\rangle.
\end{align*}
We proceed by calculating
\begin{align*}
\langle\psi,e_i\cdot\D^2\psi\rangle=&\frac{1}{3}\langle\psi,e_i\cdot\D\big(R^N(\psi,\psi)\psi\big)\rangle \\
=&\frac{1}{3}\langle\psi,e_i\cdot\big(\tilde\nabla (R^N(\psi,\psi))\big)\cdot\psi\rangle+\underbrace{\frac{1}{3}\langle\psi,e_i\cdot R^N(\psi,\psi)\D\psi\rangle}_{=-\langle\D\psi,e_i\cdot\D\psi\rangle=0},
\end{align*}
where we have used that \(\psi\) is a solution of \eqref{euler-lagrange-psi} twice. The first term on the right hand side can further be manipulated as
\begin{align*}
\langle\psi,e_i\cdot\big(\tilde\nabla (R^N(\psi,\psi))\big)\cdot\psi\rangle&=-\langle\psi,\tilde\nabla_{e_i}\big(R^N(\psi,\psi)\big)\psi\rangle \\
&=-\big\langle d\phi(e_i),\langle(\nabla R^N)^\sharp(\psi,\psi)\psi,\psi\rangle\big\rangle
-2\langle\tilde\nabla_{e_i}\psi,R^N(\psi,\psi)\psi\rangle.
\end{align*}
In addition, we find
\begin{align*}
\langle\psi,e_i\cdot e_r\cdot e_s\cdot R^N(d\phi(e_r),d\phi(e_s))\psi\rangle
=&R_{\alpha\beta\gamma\delta}\langle\psi^\alpha,e_i\cdot e_r\cdot e_s\cdot\psi^\beta\rangle_{\Sigma M}\frac{\partial\phi^\gamma}{\partial x^r}\frac{\partial\phi^\delta}{\partial x^s}\\
=&R_{\alpha\beta\gamma\delta}\langle e_i\cdot e_r\cdot e_s\cdot\psi^\alpha,\psi^\beta\rangle_{\Sigma M}\frac{\partial\phi^\gamma}{\partial x^r}\frac{\partial\phi^\delta}{\partial x^s}.
\end{align*}
A careful inspection of this term reveals that it is both real and imaginary and thus has to vanish except in the cases \(i=r\) or \(i=s\).
Consequently, we find
\begin{align*}
\langle\psi,e_i\cdot e_r\cdot e_s\cdot R^N(d\phi(e_r),d\phi(e_s))\psi\rangle=-2\langle\psi,e_r\cdot R^N(d\phi(e_i),d\phi(e_r))\psi\rangle.
\end{align*}
Combining the previous equations we find
\begin{align}
\label{identity-d}
\langle\psi,e_i\cdot\tilde\Delta\psi\rangle=&-\langle\psi,e_r\cdot R^N(d\phi(e_i),d\phi(e_r))\psi\rangle \\
\nonumber&+\frac{1}{3}\big\langle d\phi(e_i),\langle(\nabla R^N)^\sharp(\psi,\psi)\psi,\psi\rangle\big\rangle 
+\frac{2}{3}\langle\tilde\nabla_{e_i}\psi,R^N(\psi,\psi)\psi\rangle.
\end{align}
Putting together \eqref{identity-a}, \eqref{identity-b}, \eqref{identity-c} and \eqref{identity-d}
then yields the claim.
\end{proof}

For a Dirac-harmonic map with curvature term the trace of \eqref{stress-energy-dhc} can easily be computed and gives
\[
g^{ij}S_{ij}=(2-n)(|d\phi|^2+\frac{1}{6}\langle R^N(\psi,\psi)\psi,\psi\rangle).
\]
Hence, we will consider the following energy
\[
e_c(\phi,\psi):=|d\phi|^2+\frac{1}{6}\langle R^N(\psi,\psi)\psi,\psi\rangle
\]
and study its monotonicity.
Note that we need to make a case distinction as in Remark \ref{rem-monotone-quantity-dhc}
in order to obtain the positivity of \(e_c(\phi,\psi)\).

\begin{Prop}[Monotonicity formula in \(\mathbb{R}^n\)]
Let \((\phi,\psi)\) be a smooth solution of \eqref{euler-lagrange-phi}, \eqref{euler-lagrange-psi} for \(M=\mathbb{R}^n\). 
Let \(B_R(x_0)\) be a geodesic ball around the point \(x_0\in M\) and \(0<R_1<R_2\leq R\).
Then the following following monotonicity formula holds
\begin{align}
\label{dhc-mono-rn}
R_1^{2-n}\int_{B_{R_1}}e_c(\phi,\psi)dx
=&R_2^{2-n}\int_{B_{R_2}}e_c(\phi,\psi)dx \\
\nonumber&-\int_{R_1}^{R_2}\big(r^{2-n}\int_{\partial B_r}(2\big|\frac{\partial\phi}{\partial r}\big|^2
+\langle\psi,\partial_r\cdot\tilde{\nabla}_{\partial_r}\psi\rangle)d\sigma\big)dr.
\end{align}
\end{Prop}

\begin{proof}
For \(M=\mathbb{R}^n\) we choose the conformal vector field \(X=r\frac{\partial}{\partial r}\) with \(r=|x|\).
In this case we have \(\operatorname{div}(X)=n\), thus we obtain
\begin{align*}
r\int_{\partial B_r(x_0)}(2\big|\frac{\partial\phi}{\partial r}\big|^2
+\langle\psi,\partial_r\cdot\tilde{\nabla}_{\partial_r}\psi\rangle
-e_c(\phi,\psi))d\sigma
=(2-n)\int_{B_r(x_0)}(|d\phi|^2+\frac{1}{6}\langle R^N(\psi,\psi)\psi,\psi\rangle)dx,
\end{align*}
where we used \eqref{conformal-vf} and \eqref{gauss-tensor-formula}.
This can be rewritten as
\[
(2-n)\int_{B_r(x_0)}e_c(\phi,\psi)+r\int_{\partial B_r(x_0)}e_c(\phi,\psi)
=r\int_{\partial B_r(x_0)}(2\big|\frac{\partial\phi}{\partial r}\big|^2
+\langle\psi,\partial_r\cdot\tilde{\nabla}_{\partial_r}\psi\rangle)dx
\]
and by applying the coarea formula we find
\begin{align*}
\frac{d}{dr}\big(r^{2-n}\int_{B_r}e_c(\phi,\psi)dx\big)
=r^{2-n}\int_{\partial B_r}(2\big|\frac{\partial\phi}{\partial r}\big|^2
+\langle\psi,\partial_r\cdot\tilde{\nabla}_{\partial_r}\psi\rangle)d\sigma.
\end{align*}
The result then follows by integration with respect to \(r\).
\end{proof}	

\begin{Bem}
The last statement also holds if \((\phi,\psi)\) is a weak Dirac-harmonic map with curvature term, that is
\((\phi,\psi)\in\chi(M,N)\) for \(M=\mathbb{R}^n\).
It this case we can require higher integrability assumptions on \(\psi\) as in \cite{MR2544729}, Proposition 4.5
to get the following result:
Let the pair \((\phi,\psi)\) be a weak Dirac-harmonic map with curvature term in some domain \(D\subset\mathbb{R}^n\).
In addition, suppose that \(\nabla\psi\in L^p(D)\) for \(\frac{2n}{3}<p\leq n\), then
\[
R_1^{2-n}\int_{B_{R_1}}e_c(\phi,\psi)dx
\leq R_2^{2-n}\int_{B_{R_2}}e_c(\phi,\psi)dx
+C_0R_2^{3-\frac{2n}{p}}.
\]
Here, the constant \(C_0\) only depends on \(|\nabla\psi|_{L^p(D)}\).

A possible application of this monotonicity formula for stationary Dirac-harmonic maps with curvature term is to calculate the Hausdorff dimension of their singular set.
For Dirac-harmonic maps this has been carried out in \cite{MR2544729}, Proposition 4.5 and was recently extended 
to Dirac-harmonic maps with curvature term in \cite{MR4018319} and 
furthermore to Dirac-harmonic maps with curvature term coupled to a gravitino in \cite{MR3798022}.
\end{Bem}

To derive a monotonicity formula on a Riemannian manifold we again fix a point \(x_0\in M\) and consider a ball with geodesic
radius \(r=d(x_0,\cdot)\) around that point, where \(d\) denotes the Riemannian distance function.

\begin{Lem}
Let \((\phi,\psi)\) be a smooth solution of the system \eqref{euler-lagrange-phi}, \eqref{euler-lagrange-psi}.
Then the following formula holds
\begin{align}
\label{pre-monotonicity-formula-dhc}
-\Omega\int_{B_r}e_c(\phi,\psi)dx+r\int_{\partial B_r}e_c(\phi,\psi)d\sigma
=&r\int_{\partial B_r}\big(2\big|\frac{\partial\phi}{\partial r}\big|^2+\langle\psi,\partial_r\cdot\tilde{\nabla}_{\partial_r}\psi\rangle\big)d\sigma \\
\nonumber & -\int_{B_r}\hess(r^2)(e_i,e_i)\langle\psi,e_i\cdot\tilde{\nabla}_{e_i}\psi\rangle dx \\
\nonumber & -2\int_{B_r}\hess(r^2)(e_i,e_i)\langle d\phi(e_i),d\phi(e_i)\rangle\big)dx,
\end{align}
where \(\Omega:=\tr\hess(r^2)\).
\end{Lem}

\begin{proof}
We apply \eqref{gauss-tensor-formula} using \eqref{stress-energy-dhc}, which yields
\begin{align*}
&r\int_{\partial B_r}\big(2\big|\frac{\partial\phi}{\partial r}\big|^2+\langle\psi,\partial_r\cdot\tilde{\nabla}_{\partial_r}\psi\rangle\big)d\sigma
-r\int_{\partial B_r}e_c(\phi,\psi)d\sigma= 
-\int_{B_r}\tr\hess(r^2)e_c(\phi,\psi)dx \\
&+\int_{B_r}\hess(r^2)(e_i,e_j)\big(\frac{1}{2}\langle\psi,e_i\cdot\tilde{\nabla}_{e_j}\psi+e_j\cdot\tilde{\nabla}_{e_i}\psi\rangle
+2\langle d\phi(e_i),d\phi(e_j)\rangle\big)dx.
\end{align*}
Diagonalizing the Hessian of the Riemann distance function then yields the claim.
\end{proof}

Again, the presence of the Dirac-Term on the right hand side of \eqref{pre-monotonicity-formula-dhc} is 
an obstacle to a monotonicity formula.
We can try to improve the result if we assume that the solution \(\psi\) of \eqref{euler-lagrange-psi} has some additional structure.

\begin{Dfn}
We call \(\psi\in\Gamma(\Sigma M\otimes\phi^\ast TN)\) a \emph{vector twistor spinor} if it satisfies
\begin{equation}
\tilde{\nabla}_X\psi+\frac{1}{n}X\cdot\D\psi=0
\end{equation}
for all vector fields \(X\).
\end{Dfn}

\begin{Bem}
If we assume that \(\psi\) is both a vector twistor spinor and a solution of \eqref{euler-lagrange-psi} we find
\[
\tilde{\nabla}_X\psi=-\frac{1}{3n}R^N(\psi,\psi)X\cdot\psi,
\]
for all vector fields \(X\). Moreover, a direct calculation yields
\[
\partial_X\frac{1}{2}|\psi|^2=\langle\tilde{\nabla}_X\psi,\psi\rangle=-\frac{1}{3n}\langle R^N(\psi,\psi)X\cdot\psi,\psi\rangle
=-\frac{1}{3n}R_{\alpha\beta\gamma\delta}\langle\psi^\alpha,\psi^\delta\rangle\langle\psi^\beta,X\cdot \psi^\gamma\rangle.
\]
On the other hand we find
\[
\overline{R_{\alpha\beta\gamma\delta}\langle\psi^\alpha,\psi^\delta\rangle\langle\psi^\beta,X\cdot \psi^\gamma\rangle}=R_{\alpha\beta\gamma\delta}\langle\psi^\delta,\psi^\alpha\rangle\langle X\cdot \psi^\gamma,\psi^\beta\rangle
=-R_{\alpha\beta\gamma\delta}\langle\psi^\alpha,\psi^\delta\rangle\langle\psi^\beta,X\cdot \psi^\gamma\rangle.
\]
Consequently the above expression is both purely imaginary and also purely real and thus has to vanish,
meaning that \(|\psi|^2\) has constant norm. Thus, this approach does not lead to an interesting monotonicity formula. 
\end{Bem}

Only the last term on the right hand side of \eqref{pre-monotonicity-formula-dhc} has a definite sign
and we can estimate it as follows
\begin{align*}
0\leq\hess(r^2)(e_i,e_i)\langle d\phi(e_i),d\phi(e_i)\rangle\leq\omega_{max}|d\phi|^2,
\end{align*}
where \(\omega_{max}\) denotes the largest eigenvalue of \(\hess(r^2)\).

Without loss of generality we assume that \(\omega_1=\omega_{\max}\) and rewrite
\begin{align*}
\hess(r^2)(e_i,e_i)\langle\psi,e_i\cdot\tilde{\nabla}_{e_i}\psi\rangle
=&\omega_{max}\langle\psi,\D\psi\rangle +\sum_{j=2}^n\langle\psi,e_j\cdot\tilde{\nabla}_{e_j}\psi\rangle(\omega_j-\omega_{max}).
\end{align*}
Using \eqref{euler-lagrange-psi} this gives us the following inequality
\begin{align}
\label{omegamax}
-2\hess(r^2)(e_i,e_i)\langle d\phi(e_i),d\phi(e_i)\rangle-\hess(r^2)(e_i,e_i)\langle\psi,e_i\cdot\tilde{\nabla}_{e_i}\psi\rangle\\
\nonumber\geq-2\omega_{max} e_c(\phi,\psi)-\sum_{j=2}^n\langle\psi,e_j\cdot\tilde{\nabla}_{e_j}\psi\rangle(\omega_j-\omega_{max}).
\end{align}

\begin{Prop}
Let \((\phi,\psi)\) be a smooth solution of the system \eqref{euler-lagrange-phi}, \eqref{euler-lagrange-psi}.
Then for all \(0<R_1<R_2\leq R\), \(R\in(0,i_M)\) the following monotonicity type formula holds
\begin{align}
\label{prea-monotonicity-formula-dhc}
R_1^{2\omega_{max}-\Omega}&\int_{B_{R_1}}\big(e_c(\phi,\psi) 
-\langle\psi,\partial_r\cdot\tilde{\nabla}_{\partial_r}\psi\rangle\big)dx \\
\nonumber\leq
&R_2^{2\omega_{max}-\Omega}\int_{B_{R_2}}\big(e_c(\phi,\psi)
-\langle\psi,\partial_r\cdot\tilde{\nabla}_{\partial_r}\psi\rangle\big)dx \\
\nonumber&+(2\omega_{max}-\Omega)\int_{R_1}^{R_2}\big(r^{2\omega_{max}-\Omega-1}\int_{B_r}
\langle\psi,\partial_r\cdot\tilde{\nabla}_{\partial_r}\psi\rangle dx\big)dr\\
\nonumber&+\int_{R_1}^{R_2}\big(r^{2\omega_{max}-\Omega-1}
\sum_{j=2}^n\int_{B_r}\langle\psi,e_j\cdot\tilde{\nabla}_{e_j}\psi\rangle(\omega_j-\omega_{max})dx\big)dr.
\end{align}

\end{Prop}
\begin{proof}
Combining \eqref{pre-monotonicity-formula-dhc} and \eqref{omegamax} we find
\begin{align*}
(2\omega_{max}-\Omega)\int_{B_r}e_c(\phi,\psi)dx+r\int_{\partial B_r}e_c(\phi,\psi)d\sigma
\geq&r\int_{\partial B_r}\big(2\big|\frac{\partial\phi}{\partial r}\big|^2+\langle\psi,\partial_r\cdot\tilde{\nabla}_{\partial_r}\psi\rangle\big)d\sigma \\
&-\sum_{j=2}^n\int_{B_r}\langle\psi,e_j\cdot\tilde{\nabla}_{e_j}\psi\rangle(\omega_j-\omega_{max})dx.
\end{align*}
Making use of the coarea formula this can be rewritten as
\begin{align*}
\frac{d}{dr}r^{2\omega_{max}-\Omega}\int_{B_r}e_c(\phi,\psi)dx\geq &
r^{2\omega_{max}-\Omega}\int_{\partial B_r}\big(2\big|\frac{\partial\phi}{\partial r}\big|^2+\langle\psi,\partial_r\cdot\tilde{\nabla}_{\partial_r}\psi\rangle\big)d\sigma \\
&-r^{2\omega_{max}-\Omega-1}\sum_{j=2}^n\int_{B_r}\langle\psi,e_j\cdot\tilde{\nabla}_{e_j}\psi\rangle(\omega_j-\omega_{max})dx \\
\geq &r^{2\omega_{max}-\Omega}\int_{\partial B_r}\langle\psi,\partial_r\cdot\tilde{\nabla}_{\partial_r}\psi\rangle d\sigma \\
&-r^{2\omega_{max}-\Omega-1}\sum_{j=2}^n\int_{B_r}\langle\psi,e_j\cdot\tilde{\nabla}_{e_j}\psi\rangle(\omega_j-\omega_{max})dx.
\end{align*}
Integrating with respect to \(r\) and using integration by parts
\begin{align*}
\int_{R_1}^{R_2}&\big(r^{2\omega_{max}-\Omega}\int_{\partial B_r}\langle\psi,\partial_r\cdot\tilde{\nabla}_{\partial_r}\psi\rangle d\sigma\big)dr \\
=&\int_{R_1}^{R_2}\big(r^{2\omega_{max}-\Omega}\frac{d}{dr}\int_{B_r}\langle\psi,\partial_r\cdot\tilde{\nabla}_{\partial_r}\psi\rangle d\sigma\big)dr \\
=&R_2^{2\omega_{max}-\Omega}\int_{B_{R_2}}\langle\psi,\partial_r\cdot\tilde{\nabla}_{\partial_r}\psi\rangle dx 
-R_1^{2\omega_{max}-\Omega}\int_{B_{R_1}}\langle\psi,\partial_r\cdot\tilde{\nabla}_{\partial_r}\psi\rangle dx \\
&-(2\omega_{max}-\Omega)\int_{R_1}^{R_2}r^{2\omega_{max}-\Omega-1}\int_{B_r}\langle\psi,\partial_r\cdot\tilde{\nabla}_{\partial_r}\psi\rangle dx
\end{align*}
completes the proof.
\end{proof}

\begin{Bem}
In the case of \(M=\mathbb{R}^n\) we have \(\omega_i=1,i=1\ldots,n\) and \(\Omega=n\).
In this case we have equality in \eqref{prea-monotonicity-formula-dhc}
and \eqref{prea-monotonicity-formula-dhc} reduces to \eqref{dhc-mono-rn}.
\end{Bem}

\begin{Bem}
Again, it seems very difficult to obtain a Liouville Theorem from the monotonicity formula \eqref{prea-monotonicity-formula-dhc}
without posing a lot of restrictions on the solution.
\end{Bem}

\subsection{A Liouville Theorem for a domain with positive Ricci curvature}
In this section we derive a vanishing theorem for Dirac-harmonic maps with
curvature term under an energy and curvature assumption, similar to Theorem \ref{theorem-liouville-complete-soler}.
To this end we set 
\[
e(\phi,\psi):=\frac{1}{2}(|d\phi|^2+|\psi|^4).
\]
\begin{Lem}[Bochner formulas]
Let \((\phi,\psi)\) be a smooth solution of the system \eqref{euler-lagrange-phi}, \eqref{euler-lagrange-psi}.
Then the following Bochner formulas hold
\begin{align}
\Delta\frac{1}{2}|\psi|^4=&\big|d|\psi|^2\big|^2+2|\psi|^2|\tilde{\nabla}\psi|^2
+\frac{2}{9}|\psi|^2|R^N(\psi,\psi)\psi|^2+\frac{R}{2}|\psi|^4 \\
\nonumber&+\frac{1}{2}|\psi|^2\langle e_i\cdot e_j\cdot R^N(d\phi(e_i),d\phi(e_j))\psi,\psi\rangle, \\
\Delta\frac{1}{2}|d\phi|^2=&|\nabla d\phi|^2+\langle d\phi(\text{Ric}^M(e_i)),d\phi(e_i)\rangle
-\langle R^N(d\phi(e_i),d\phi(e_j))d\phi(e_j),d\phi(e_i)\rangle \\
\nonumber&+\frac{1}{2}\langle(\nabla_{d\phi(e_j)}R^N)(\psi,e_i\cdot\psi)d\phi(e_i),d\phi(e_j)\rangle
+\langle R^N(\psi,e_i\cdot\tilde{\nabla}_{e_j}\psi)d\phi(e_i),d\phi(e_j)\rangle \\
\nonumber&+\frac{1}{2}\langle R^N(\psi,e_i\cdot\psi)\nabla_{e_j}d\phi(e_i),d\phi(e_j)\rangle
+\frac{1}{12}\langle\langle(\nabla_{d\phi(e_i)}(\nabla R^N)^\sharp)(\psi,\psi)\psi,\psi\rangle,d\phi(e_i)\rangle \\
\nonumber&+\frac{1}{3}\langle\langle(\nabla R^N)^\sharp(\tilde\nabla_{e_i}\psi,\psi)\psi,\psi\rangle,d\phi(e_i)\rangle,
\end{align}
where \(e_i,i=1\ldots,n\) is an orthonormal basis of \(TM\).
\end{Lem}

\begin{proof}
We choose a local orthonormal basis of \(TM\) such that \(\nabla_{e_i}e_j=0,i,j=1,\ldots,n\) at the considered point.
The fist equation follows by a direct calculation using the Weitzenböck formula for the twisted Dirac-operator \(\D\), that is
\begin{equation*}
\D^2\psi=-\tilde{\Delta}\psi+\frac{1}{4}R\psi +\frac{1}{2}e_i\cdot e_j\cdot R^N(d\phi(e_i),d\phi(e_j))\psi,
\end{equation*}
where \(\tilde{\Delta}\) denotes the connection Laplacian on the vector bundle \(\Sigma M\otimes\phi^\ast TN\).
To obtain the second equation we recall the following Bochner formula
for a map \(\phi\colon M\to N\)
\[
\Delta\frac{1}{2}|d\phi|^2=|\nabla d\phi|^2+\langle d\phi(\text{Ric}^M(e_i)),d\phi(e_i)\rangle
-\langle R^N(d\phi(e_i),d\phi(e_j))d\phi(e_j),d\phi(e_i)\rangle+\langle\nabla\tau(\phi),d\phi\rangle.
\]
Moreover, by a direct calculation we obtain
\begin{align*}
\tilde\nabla_{e_j}\big(\frac{1}{2}R^N(\psi,e_i\cdot\psi)d\phi(e_i)\big)=&
\frac{1}{2}(\nabla_{d\phi(e_j)}R^N)(\psi,e_i\cdot\psi)d\phi(e_i)+R^N(\psi,e_i\cdot\tilde{\nabla}_{e_j}\psi)d\phi(e_i) \\
&+\frac{1}{2}R^N(\psi,e_i\cdot\psi)\nabla_{e_j}d\phi(e_i), \\
\tilde\nabla_{e_j}\big(\frac{1}{12}\langle(\nabla R^N)^\sharp(\psi,\psi)\psi,\psi\rangle\big)
=&\frac{1}{12}\langle(\nabla_{d\phi(e_j)}(\nabla R^N)^\sharp)(\psi,\psi)\psi,\psi\rangle
+\frac{1}{3}\langle(\nabla R^N)^\sharp(\tilde\nabla_{e_j}\psi,\psi)\psi,\psi\rangle,
\end{align*}
which concludes the proof.
\end{proof}

\begin{Cor}
Let \((\phi,\psi)\) be a smooth solution of the system \eqref{euler-lagrange-phi}, \eqref{euler-lagrange-psi}.
Then the following estimate holds:
\begin{equation}
\label{bochner-energy}
\Delta e(\phi,\psi)\geq c_1(|\nabla d\phi|^2+\big|d|\psi|^2\big|^2)-c_2e(\phi,\psi)-c_3(e(\phi,\psi))^2,
\end{equation}
where \(c_i,i=1,2,3\) are positive constants that depend only on the geometry of \(M\) and \(N\).
\end{Cor}

\begin{proof}
Making use of the Bochner formulas we find
\begin{align*}
\Delta e(\phi,\psi)\geq& |\nabla d\phi|^2+\kappa_M|d\phi|^2+\kappa_N|d\phi|^4
-\frac{|\nabla R^N|_{L^\infty}\sqrt{n}}{2}|\psi|^2|d\phi|^3-|R^N|_{L^\infty}\sqrt{n}|\psi||\tilde{\nabla}\psi||d\phi|^2 \\
&-\frac{|R^N|_{L^\infty}\sqrt{n}}{2}|\psi|^2|\nabla d\phi||d\phi|
-\frac{|\nabla^2R^N|_{L^\infty}}{12}|\psi|^4|d\phi|^2
-\frac{|\nabla R^N|_{L^\infty}}{3}|\tilde{\nabla}\psi||\psi|^3|d\phi| \\
&
+\big|d|\psi|^2\big|^2+2|\psi|^2|\tilde{\nabla}\psi|^2
+\frac{2}{9}|\psi|^2|R^N(\psi,\psi)\psi|^2+\frac{R}{2}|\psi|^4 
-\frac{n|R^N|_{L^\infty}}{2}|\psi|^4|d\phi|^2, 
\end{align*}
where \(\kappa_M\) denotes a lower bound for the Ricci curvature of \(M\) and
\(\kappa_N\) an upper bound for the sectional curvature of \(N\).
By application of Young's inequality we find
\begin{align}
\Delta e(\phi,\psi)
\label{inequality-dhc-deltas}\geq & (1-\delta_1)|\nabla d\phi|^2+|\psi|^2|\tilde{\nabla}\psi|^2(2-\delta_2-\delta_3)+\big|d|\psi|^2\big|^2\\
\nonumber&+\frac{2}{9}|\psi|^2|R^N(\psi,\psi)\psi|^2+\frac{R}{2}|\psi|^4+\kappa_M|d\phi|^2 \\
\nonumber&-|d\phi|^4 (-\kappa_N+\frac{1}{\delta_2}\frac{n}{4}|R^N|^2_{L^\infty}+\delta_4) \\
\nonumber&-|\psi|^4|d\phi|^2\big(\frac{n}{8\delta_1}|R^N|^2_{L^\infty}+\frac{1}{\delta_3}\frac{|\nabla R^N|^2_{L^\infty}}{36}
+\frac{|\nabla^2R^N|_{L^\infty}}{12}+\frac{n|R^N|_{L^\infty}}{2}+\frac{n}{8\delta_4}|\nabla R^N|^2_{L^\infty}\big) 
\end{align}
with positive constants \(\delta_i,i=1,\ldots 4\).
The statement then follows by applying Young's inequality again.
\end{proof}

\begin{Bem}
\begin{enumerate}
 \item The analytic structure of \eqref{bochner-energy} is the same as in the case of harmonic maps.
 \item If we want to derive a Liouville Theorem from \eqref{bochner-energy} making only assumptions on the geometry of \(M\)
  and \(N\) we would require that both \(c_2\leq 0\) and \(c_3\leq 0\). However, it can easily be checked that we cannot achieve such an estimate
   since the curvature tensor of \(N\) appears on the right hand side of the system \eqref{euler-lagrange-phi} and \eqref{euler-lagrange-psi}.
\end{enumerate}
\end{Bem}

However, we can give a Liouville theorem under similar assumptions as in Theorem \ref{theorem-liouville-complete-soler}.
A similar Theorem for Dirac-harmonic maps was obtained in \cite{MR3144359}, Theorem 4.
\begin{Satz}
\label{theorem-liouville-dhc}
Let \((M,g)\) be a complete noncompact Riemannian spin manifold and \((N,h)\) be a Riemannian manifold
with 	nonpositive curvature. Suppose that \((\phi,\psi)\) is a Dirac-harmonic map with curvature term 
with finite energy \(e(\phi,\psi)\).
If
\begin{equation}
\label{assumption-liouville-complete-dhc}
\operatorname{Ric}_M\geq (c_1|\psi|^4+c_2|d\phi|^2)g,
\end{equation}
with the constants
\begin{align*}
c_1=&\frac{n}{2}|R^N|_{L^\infty}+\frac{n}{16}|R^N|^2_{L^\infty}+\big(\frac{1}{36}+\frac{n}{8}\big)|\nabla R^N|^2_{L^\infty}
+\frac{|\nabla^2R^N|_{L^\infty}}{12},\\
c_2=&\frac{n}{4}|R^N|^2_{L^\infty}+1
\end{align*}
then \(\phi\) maps to a point and \(\psi\) vanishes identically.
\end{Satz}

\begin{proof}
First of all we note that 
\begin{align}
\label{liouville-dhc-a}
|d e(\phi,\psi)|^2=|\frac{1}{2}d(|d\phi|^2+|\psi|^4)|^2\leq&
(|d\phi|^2|\nabla d\phi|^2+|\psi|^4\big|d|\psi|^2\big|^2+2|d\phi||\nabla d\phi||\psi|^2\big|d|\psi|^2\big|) \\
\nonumber\leq &2e(\phi,\psi)(|\nabla d\phi|^2+\big|d|\psi|^2\big|^2).
\end{align}
If we put \(\delta_1=\frac{1}{2},\delta_2=\delta_3=\delta_4=1\) in \eqref{inequality-dhc-deltas} we find
\begin{align*}
\Delta e(\phi,\psi)\geq&\frac{1}{2}(|\nabla d\phi|^2+\big|d|\psi|^2\big|^2) \\
&+|d\phi|^2\bigg(\kappa_M-|\psi|^4\big(\frac{n}{2}|R^N|_{L^\infty}+\frac{n}{16}|R^N|^2_{L^\infty}+\big(\frac{1}{36}+\frac{n}{8}\big)|\nabla R^N|^2_{L^\infty}
+\frac{|\nabla^2R^N|_{L^\infty}}{12}\big) \\
&+|d\phi|^2\big(1+\frac{n}{4}|R^N|^2_{L^\infty}\big)\bigg)
\end{align*}
Making use of the assumption \eqref{assumption-liouville-complete-dhc} this yields
\begin{equation}
\label{liouville-dhc-b}
\Delta e(\phi,\psi)\geq \delta(|\nabla d\phi|^2+\big|d|\psi|^2\big|^2)
\end{equation}
for a positive constant \(\delta\). We fix a positive number \(\epsilon>0\) such that
\begin{align*}
\Delta\sqrt{\delta e(\phi,\psi)+\epsilon}=&\frac{\delta\Delta e(\phi,\psi)}{2\sqrt{e(\phi,\psi)+\epsilon}}
-\frac{1}{4}\frac{\delta^2|\nabla e(\phi,\psi)|^2}{(e(\phi,\psi)+\epsilon)^\frac{3}{2}} \\
\geq&\delta^2\frac{|\nabla d\phi|^2+\big|d|\psi|^2\big|^2}{2\sqrt{e(\phi,\psi)+\epsilon}}\big(1-\frac{e(\phi,\psi)}{e(\phi,\psi)+\epsilon}\big) 
\geq 0,
\end{align*}
where we used \eqref{liouville-dhc-a} and \eqref{liouville-dhc-b}.
The rest of the proof is identical to the proof of Theorem \ref{theorem-liouville-complete-soler}.
\end{proof}

\par\medskip
\textbf{Acknowledgements:}
The author gratefully acknowledges the support of the Austrian Science Fund (FWF) 
through the project P30749-N35 ``Geometric variational problems from string theory''.

\bibliographystyle{plain}
\bibliography{mybib}
\end{document}